\documentclass[11pt]{article}

\usepackage[margin=1in]{geometry}
\usepackage[T1]{fontenc}
\usepackage[utf8]{inputenc}
\usepackage{lmodern}
\usepackage{microtype}
\usepackage{amsmath}
\usepackage{amssymb}
\usepackage{amsthm}
\usepackage{mathtools}
\usepackage{graphicx}
\usepackage{booktabs}
\usepackage{multirow}
\usepackage{array}
\usepackage{makecell}
\usepackage{float}
\usepackage{caption}
\usepackage{subcaption}
\usepackage{algorithm}
\usepackage{algpseudocode}
\usepackage{xcolor}
\usepackage{url}
\usepackage[colorlinks=true,
            linkcolor=blue,
            citecolor=blue,
            urlcolor=blue]{hyperref}
\usepackage[nameinlink,capitalize]{cleveref}
\usepackage{listings}
\usepackage{titlesec}
\usepackage{enumitem}
\usepackage[section]{placeins}
\usepackage{needspace}
\usepackage{etoolbox}

\titlespacing*{\section}{0pt}{2.2ex plus .5ex minus .2ex}{1.1ex plus .2ex}
\titlespacing*{\subsection}{0pt}{1.8ex plus .4ex minus .2ex}{.8ex plus .2ex}
\titlespacing*{\paragraph}{0pt}{1.1ex plus .3ex minus .2ex}{.8em}
\setlist[itemize]{topsep=4pt,itemsep=3pt,parsep=0pt,partopsep=0pt,leftmargin=*}
\setlist[enumerate]{topsep=4pt,itemsep=5pt,parsep=0pt,partopsep=0pt,leftmargin=*}
\AtBeginDocument{%
  \setlength{\abovedisplayskip}{8pt plus 2pt minus 2pt}%
  \setlength{\belowdisplayskip}{8pt plus 2pt minus 2pt}%
  \setlength{\abovedisplayshortskip}{3pt plus 2pt}%
  \setlength{\belowdisplayshortskip}{5pt plus 2pt minus 2pt}%
}
\AtBeginEnvironment{algorithm}{\small}
\makeatletter
\def\theHALG@line{\thealgorithm.\arabic{ALG@line}}
\makeatother
\definecolor{codegray}{gray}{0.95}
\definecolor{codegreen}{rgb}{0.0,0.45,0.0}
\definecolor{codepurple}{rgb}{0.45,0.0,0.55}

\definecolor{revcolor}{RGB}{180,0,0} 
\newcommand{\revise}[1]{{\color{black}#1}}

\lstdefinestyle{pythonstyle}{
    language=Python,
    aboveskip=8pt,
    belowskip=8pt,
    backgroundcolor=\color{codegray},
    basicstyle=\ttfamily\small,
    keywordstyle=\color{blue},
    commentstyle=\color{codegreen},
    stringstyle=\color{codepurple},
    numbers=left,
    numberstyle=\tiny,
    stepnumber=1,
    numbersep=8pt,
    showstringspaces=false,
    breaklines=true,
    columns=fullflexible,
    keepspaces=true,
    frame=single,
    tabsize=4,
    captionpos=b
}

\usepackage{tikz}
\usepackage{pgfplots}
\pgfplotsset{compat=1.18}

\newtheorem{proposition}{Proposition}

\theoremstyle{definition}

\DeclareMathOperator{\rank}{rank}
\DeclareMathOperator{\supp}{supp}

\DeclareMathOperator{\xc}{xc}
\DeclareMathOperator{\bin}{bin}
\newcommand{\Sectionref}[1]{Section~\ref{#1}}
\newcommand{\Figureref}[1]{Figure~\ref{#1}}
\newcommand{\Tableref}[1]{Table~\ref{#1}}
\newcommand{\Appendixref}[1]{Appendix~\ref{#1}}

\title{Computing Lower Bounds on the Nonnegative Rank via \\ Non-Convex Optimization Solvers}

\author{Timothy Baeckelant, Arnaud Vandaele, Nicolas Gillis\thanks{Department of Mathematics and Operational Research, University of Mons, Mons, Belgium. We acknowledge the support by the European Union (ERC consolidator, eLinoR, no 101085607). \\ Emails: firstname.lastname@umons.ac.be.}}

\date{}

\begin{document}
\maketitle

\begin{abstract}
The \textit{nonnegative rank} of a nonnegative matrix $X$ is the smallest number of nonnegative rank-one factors that sum to $X$.
Since computing the nonnegative rank is NP-hard, it is common to circumvent this issue by computing lower and upper bounds.
In this paper, we propose non-convex formulations and practical implementations for four important lower bounds for the nonnegative rank, namely the \textit{fooling set} bound (FSB), the \textit{rectangle covering} bound (RCB), the \textit{hyperplane separation} bound (HSB), and the \textit{self-scaled} bound (SSB).
In particular, to the best of our knowledge, our algorithm for computing the SSB is the first available in the literature.
We compute the SSB for several matrices for which its value had not previously been reported.
In some cases, the SSB provides the best known upper bound, while the other three bounds do not.
Moreover, on canonical benchmarks, we show that our non-convex approaches are often competitive with existing methods.
The paper also provides a consolidated reference for the current state of several classical lower bounds on a large number of benchmark matrices.
\end{abstract}


\section{Introduction}\label{sec:intro}

Given a nonnegative matrix $X\in\mathbb{R}_+^{m\times n}$ and a target rank $r\in\mathbb{N}$, nonnegative matrix factorization (NMF) seeks $W\in\mathbb{R}_+^{m\times r}$ and $H\in\mathbb{R}_+^{r\times n}$ such that $X\approx WH$. Approximating a data matrix $X$ by a low-rank nonnegative product often leads to meaningful part-based representations with numerous applications in machine learning and data analysis~\cite{cichocki2009nonnegative, gillis2020nonnegative}. When $X=WH$ holds exactly, we say that $X$ admits an Exact NMF of size $r$, and the smallest such $r$ is called the nonnegative rank of $X$ and denoted $\rank_+(X)$.

\paragraph{Computational complexity}
Unlike the usual rank of a matrix, computing $\rank_+(X)$ is challenging: deciding whether $\rank_+(X)=\rank(X)$ is NP-hard~\cite{vavasis2010complexity} and computing $\rank_+(X)$ remains NP-hard even when $X$ is a $\{0,1\}$-matrix~\cite{shitov2017nonnegative}.
When $r$ is fixed (that is, not part of the input), algorithms based on quantifier elimination decide exact factorizations in time $(mn)^{\mathcal{O}(r^2)}$~\cite{moitra2016almost}, but to the best of our knowledge, no practical implementation of such methods exists. 
In practice, one therefore relies on bounds: upper bounds via explicit factorizations, and lower bounds that certify $\rank_+(X)\ge \ell$ for some constant $\ell$, without computing a factorization. 
This paper focuses on the latter.

\paragraph{Why lower bounds matter: extension complexity in a nutshell}
Lower bounds on $\rank_+(X)$ have direct consequences in combinatorial optimization via the extension complexity.
For a polytope $P=\{x\in\mathbb{R}^d\mid Ax\le b\}$ with $v$ vertices, denoted $\{x_j\}_{j=1}^v$, and $f$ facets, $\{A(i,:)x\le b_i\}_{i=1}^f$, its slack matrix, $S_P\in\mathbb{R}_+^{f\times v}$, is defined by $S_P(i,j)=b_i - A(i,:) x_j$.
An extended formulation for $P$ is a higher dimensional polytope that linearly projects onto $P$. 
It turns out that the minimum number of facets of any extended formulation of $P$, called its extension complexity and denoted $\xc(P)$, equals the nonnegative rank of its 
slack matrix~\cite{yannakakis1988expressing}:
\[
\xc(P)=\rank_+(S_P).
\]
Thus, any lower bound on $\rank_+(S_P)$ is a lower bound on the extension complexity of $P$. 
This link between polytope theory and the nonnegative rank underlies several landmark results on the limits of linear descriptions for combinatorial polytopes~\cite{fiorini2012linear,fiorini2015exponential,kaibel2015short,rothvoss2017matching}.

\paragraph{Contributions} 
Most lower bounds for the nonnegative rank are combinatorial, depending only on the zero–nonzero pattern of $X$; see, e.g.,~\cite{fiorini2013combinatorial, fiorini2019}. A smaller number of bounds also exploit the magnitude of the entries; see,
for example, \cite{fawzi2015lower}.
We refer the reader to~\cite[Chapter~3]{gillis2020nonnegative} for an overview of lower bounds on the nonnegative rank. 
In this paper, we investigate the use of exact non-convex optimization models to compute classical lower bounds on the nonnegative rank.
Our goal is not to introduce new bounds, but rather to show that several standard bounds can be reformulated or, in the case of iterative schemes, implemented within a non-convex optimization framework that can be handled by modern global solvers.
This direction has remained largely unexplored in practice, mainly because globally solving non-convex quadratic programs has only recently become computationally realistic\footnote{For example, a non-convex quadratic solver was introduced in Gurobi 9.0 (2019).}.
The situation has changed: modern solvers make exact non-convex approaches both practical and certifiable on problem sizes of interest.
We study this approach on four representative bounds: the fooling set bound (FSB), the rectangle covering bound (RCB), the hyperplane separation bound (HSB), and the self-scaled bound (SSB).

Our contributions are threefold.
First, we derive non-convex formulations for FSB and RCB, and a non-convex rank-one separation oracle for HSB.
Second, to the best of our knowledge, we provide the first algorithm for directly computing the SSB using the same approach as for the HSB.
This implementation allows us to evaluate the SSB on several matrices for which its value had not previously been reported.
On some of these matrices, the SSB gives a larger lower bound on the nonnegative rank than the FSB, RCB, or HSB.
Third, we benchmark these non-convex approaches against existing implementations on several matrix families to identify the situations in which they are advantageous and those in which they are not.
In addition, the paper consolidates updated computational information on these bounds for several matrix families.
To this end, \Appendixref{app:tables} provides a detailed overview of the currently best known values and bounds for the four lower bounds considered in this work, together with the corresponding information on the nonnegative rank.

\paragraph{Outline of the paper} 
The four sections after this introduction, Sections~\ref{sec:fs}--\ref{sec:ssb}, present each of the bounds: specifically, the FSB, the RCB, the HSB, and the SSB. Each section defines its corresponding bound, presents the standard computational approaches, and finally proposes our non-convex reformulation. 
\Sectionref{sec:numerics} reports and analyzes the experimental results, comparing the various approaches to compute the four lower bounds and discussing the strengths and limitations of the proposed approaches compared to the state of the art.  
\Sectionref{sec:conclusion} concludes and discusses directions for further work. 

\paragraph{Running example} To illustrate our constructions in Sections~\ref{sec:fs}--\ref{sec:ssb}, we will use the linear Euclidean distance matrices $M_n\in\mathbb{R}_+^{n\times n}$  defined entrywise by $M_n(i,j)=(i-j)^2$, for $n\in\mathbb{N}$.
More precisely, we will use $M_5$, while $M_n$ will be used in our benchmark for the numerical results presented in \Sectionref{sec:numerics}. 
These LEDMs are dense matrices, with zeros only on the diagonal, and satisfy $\rank(M_n)=3$ for all $n\ge 3$, while $\rank_+(M_n) \geq \min   \big\{r \ | \ \binom{r}{\lfloor r/2 \rfloor} \geq n \big\}$~\cite{de1981boolean}, making them a convenient structured testbed for lower-bound methods. 

\section{Fooling set bound (FSB)}\label{sec:fs}

The FSB relies on identifying positive entries that cannot be covered simultaneously by the same nonnegative rank-one factor without violating the zero pattern of the matrix. The maximum number of such mutually incompatible entries provides a lower bound on the nonnegative rank.
Let $X\in\mathbb{R}_+^{m\times n}$, and denote by $\supp(X)=\{(i,j)\mid X_{ij}>0\}$ its support.
Two positive entries of $X$ on different rows and columns, $X(i,j) > 0$ and $X(k,l) > 0$ with $i\neq k$ and $j\neq l$, are said to be pairwise independent if at least one of the two cross entries is zero, that is, if $X_{il}=0$ or $X_{kj}=0$. 
A fooling set is a set of positive entries of $X$ that are pairwise independent.

Suppose now that $X=\sum_{p=1}^r W(:,p)H(p,:)$ is an exact NMF of $X$. 
Two positive entries of $X$ at positions $(i,j)$ and $(k,l)$ are said to be covered by the same rank-one factor, say $W(:,p)H(p,:)$, if $W_{ip}H_{pj}>0$ and $W_{kp}H_{pl}>0$. If $(i,j)$ and $(k,l)$ are part of a fooling set, this would imply that both cross entries,  $X_{il}$ and $X_{kj}$, are positive since $W_{ip}H_{pl}>0$ and $W_{kp}H_{pj}>0$. Since there are no cancellations in nonnegative factorizations, this leads to a contradiction. 
Therefore, each rank-one factor can cover at most one entry of a fooling set.
The quantity $\mathrm{FSB}(X)$ denotes the maximum cardinality of a fooling set of $X$ and $\rank_+(X)\geq \mathrm{FSB}(X)$.

\subsection{Standard approach via maximum independent-set} 

The standard way to compute the FSB is to define a binary variable $x_{ij}$ for each $(i,j) \in \supp(X)$, and $x_{ij} = 1$ if $(i,j)$ is selected in the fooling set, $x_{ij} = 0$ otherwise.  
Based on a maximum independent-set formulation where each conflict is represented by a linear constraint, the maximum fooling set is obtained as the solution to the following mixed integer linear program (MILP):   
\begin{equation}\label{eq:fsb-enumFio}
\begin{aligned}
\max_{x\in\{0,1\}^{|\supp(X)|}}\quad
    & \sum_{(i,j)\in\supp(X)}x_{ij} \\
\text{such that}\quad
    & x_{ij}+x_{kl}\leq 1, \\
    & \quad\text{for all }(i,j),(k,l)\in\supp(X),
      \quad(i,j)\neq(k,l), \\
    & \quad\text{such that }(i,l),(k,j)\in\supp(X).
\end{aligned}
\end{equation}
This model uses $|\supp(X)|$ binary variables and $\mathcal{O}(|\supp(X)|^2)$ linear constraints.
It is the formulation used, for instance, in the package~\cite{fiorini2019}.

\subsection{Proposed non-convex formulation}

The standard formulation may become heavy on matrices with dense support, as it requires a quadratic number of pairwise conflict constraints in the number of non-zeros. 
Our reformulation follows a different philosophy.
Instead of introducing one binary variable per admissible entry and explicitly encoding all pairwise conflicts, we directly try to detect whether a fooling set of size $r$ exists. 
Let us define 
\[
Z=\bin(X)\in\{0,1\}^{m\times n},
\]
where $Z_{ij}=1$ if and only if $X_{ij}>0$.
We represent the $r$ selected entries through two binary assignment matrices $U\in\{0,1\}^{m\times r}$ and $V\in\{0,1\}^{r\times n}$. 
We denote by $[r] := \{1,\dots,r\}$ the set of indices from $1$ to $r$.
For each $k\in[r]$, the $k$th column of $U$ selects the row index $i_k$ of the $k$th chosen
entry, while the $k$th row of $V$ selects its column index $j_k$.
The constraints must then ensure that:
(i)~each selected position $(i_k,j_k)$ belongs to $\supp(X)$,
(ii)~no two selected entries share the same row or the same column, and 
(iii)~for every pair of selected entries, at least one of the two cross positions lies outside the support of $X$.
This leads to the following feasibility problem:
\begin{equation}\label{eq:fsb-nc}
\begin{aligned}
&\text{Find}\quad  U\in\{0,1\}^{m\times r},\quad V\in\{0,1\}^{r\times n} \quad \text{such that }\\
& \sum_{\substack{(i,j):\, Z_{ij}=1}}
\bigl(U(i,k)V(l,j)+U(i,l)V(k,j)\bigr)\le 1,
\qquad 1\leq k<l\leq r, \\
& 
\sum_{i=1}^m U(i,k)=1, \; \sum_{j=1}^n V(k,j)=1, \text{ for all } k \in [r], \\
& \sum_{k=1}^r U(i,k)\leq 1 \text{ for all } i \in [m], \; 
\sum_{k=1}^r V(k,j)\leq 1 \text{ for all } j \in [n], \\
& U(i,k)V(k,j)=0  \text{ for all } (i,j)\notin\supp(X), k \in [r].
\end{aligned}
\end{equation}
The first constraint enforces the fooling-set cross condition between any two selected entries: for the two entries encoded by indices $k$ and $l$, it ensures that at most one of the two cross positions belongs to the support of $X$.
The next four constraints are assignment constraints that ensure that each selected entry uses exactly
one row and one column, and that no row or column is used twice.
The last admissibility constraints forbid the selection of entries outside the support of $X$.
The model uses $mr+rn$ binary variables, $2r+m+n$ linear assignment constraints, $\binom{r}{2}$
bilinear cross constraints, and $r\bigl(mn-|\supp(X)|\bigr)$ bilinear admissibility constraints.
In particular, when $X$ is dense, the number of admissibility constraints remains small, in contrast to the standard formulation, whose number of conflict constraints grows quadratically with $|\supp(X)|$.
This makes the proposed model particularly attractive for dense matrices.
If the feasibility problem is satisfiable, then it returns a fooling set of size $r$. \Figureref{fig:ncfsLEDM} illustrates such a solution on the matrix $M_5$ for which there exists a fooling set of size $3$.

\begin{figure}[!htbp]
    \centering
\begin{tikzpicture}[scale=0.9, every node/.style={font=\small}]
    \def\s{0.42} 

    \begin{scope}[shift={(0,0)}]
        \node[font=\small\bfseries] at (1.05,-2.5) {$Z=\bin(M_5)$};

        \foreach \i in {0,...,4}{
            \foreach \j in {0,...,4}{
                \pgfmathtruncatemacro{\row}{\i+1}
                \pgfmathtruncatemacro{\col}{\j+1}
                \ifnum\row=\col
                    \fill[red!10] (\j*\s,-\i*\s) rectangle ++(\s,-\s);
                    \draw[red!70!black, thick] (\j*\s,-\i*\s) rectangle ++(\s,-\s);
                \else
                    \fill[gray!15] (\j*\s,-\i*\s) rectangle ++(\s,-\s);
                    \draw[black!50] (\j*\s,-\i*\s) rectangle ++(\s,-\s);
                \fi
            }
        }

        \foreach \r/\c in {1/2,2/3,3/1}{
            \pgfmathsetmacro{\x}{(\c-1)*\s}
            \pgfmathsetmacro{\y}{-(\r-1)*\s}
            \fill[blue!55] (\x,\y) rectangle ++(\s,-\s);
            \draw[blue!80!black, thick] (\x,\y) rectangle ++(\s,-\s);
        }

        \foreach \i in {1,...,5}{
            \node[left]  at (-0.08,{-((\i-1)*\s+\s/2)}) {\scriptsize \i};
            \node[above] at ({(\i-1)*\s+\s/2},0.08) {\scriptsize \i};
        }

    \end{scope}

    \begin{scope}[shift={(4.2,0)}]
        \node[font=\small\bfseries] at (0.63,-2.5) {$U$};

        \foreach \i in {0,...,4}{
            \foreach \j in {0,...,2}{
                \draw[black!50] (\j*\s,-\i*\s) rectangle ++(\s,-\s);
            }
        }

        \foreach \r/\c in {1/1,2/2,3/3}{
            \pgfmathsetmacro{\x}{(\c-1)*\s}
            \pgfmathsetmacro{\y}{-(\r-1)*\s}
            \fill[blue!55] (\x,\y) rectangle ++(\s,-\s);
            \draw[blue!80!black, thick] (\x,\y) rectangle ++(\s,-\s);
            \node at (\x+0.5*\s,\y-0.5*\s) {\scriptsize 1};
        }

        \foreach \i in {1,...,5}{
            \node[left] at (-0.08,{-((\i-1)*\s+\s/2)}) {\scriptsize \i};
        }
        \foreach \j in {1,...,3}{
            \node[above] at ({(\j-1)*\s+\s/2},0.08) {\scriptsize \j};
        }
    \end{scope}

    \begin{scope}[shift={(6.4,0)}]
        \node[font=\small\bfseries] at (1.05,-2.5) {$V$};

        \foreach \i in {0,...,2}{
            \foreach \j in {0,...,4}{
                \draw[black!50] (\j*\s,-\i*\s) rectangle ++(\s,-\s);
            }
        }

        \foreach \r/\c in {1/2,2/3,3/1}{
            \pgfmathsetmacro{\x}{(\c-1)*\s}
            \pgfmathsetmacro{\y}{-(\r-1)*\s}
            \fill[blue!55] (\x,\y) rectangle ++(\s,-\s);
            \draw[blue!80!black, thick] (\x,\y) rectangle ++(\s,-\s);
            \node at (\x+0.5*\s,\y-0.5*\s) {\scriptsize 1};
        }

        \foreach \i in {1,...,3}{
            \node[left] at (-0.08,{-((\i-1)*\s+\s/2)}) {\scriptsize \i};
        }
        \foreach \j in {1,...,5}{
            \node[above] at ({(\j-1)*\s+\s/2},0.08) {\scriptsize \j};
        }
    \end{scope}
\end{tikzpicture}
    \caption{Illustration of 
    the FSB computation using~\eqref{eq:fsb-nc}, with $\mathrm{FSB}(M_5)=3$. 
    Left: the support matrix $Z=\bin(M_5)$. The three blue entries, $(1,2)$, $(2,3)$, and $(3,1)$, form a fooling set of size $3$, while the red diagonal entries are the cross-zero condition for each pair. Right: binary assignment matrices $U$ and $V$ encoding the selected rows and columns.}
    \label{fig:ncfsLEDM}
\end{figure}

\section{Rectangle covering bound (RCB)}\label{sec:rc}

The RCB is based on the observation that the support of a nonnegative rank-one matrix forms a \textit{rectangle}. 
The RCB measures how many rectangles are needed to cover the support of the matrix.
Let $X\in\mathbb{R}_+^{m\times n}$, and let $Z=\bin(X)\in\{0,1\}^{m\times n}$.
A rectangle supported by $Z$ is a Cartesian product $I\times J\subseteq [m]\times[n]$ such that $I\times J \subseteq \supp(Z)$.
Equivalently, a rectangle can be identified with a binary
rank-one matrix $uv^\top\in\{0,1\}^{m\times n}$, where $u\in\{0,1\}^m$, $v\in\{0,1\}^n$, and $\supp(uv^\top)\subseteq \supp(Z)$.
We denote by
\[
\mathcal{R}_{01}(Z)=\bigl\{uv^\top\in\{0,1\}^{m\times n}: u\in\{0,1\}^m,\ v\in\{0,1\}^n,\ \supp(uv^\top)\subseteq \supp(Z)\bigr\}
\]
the set of all such binary rank-one rectangles.
A rectangle cover of $Z$ is a collection of rectangles in $\mathcal{R}_{01}(Z)$ whose union is exactly $\text{supp}(Z)$.
The rectangle covering number, denoted by $\mathrm{RCB}(X)$, is the minimum number of rectangles
in such a cover: 
\[
\mathrm{RCB}(X)=\min\bigl\{r:\exists U\in\{0,1\}^{m\times r},\ V\in\{0,1\}^{r\times n}\ \text{s.t.}\ \supp(Z)=\supp(UV)\bigr\}.
\] 
The RCB coincides with the Boolean rank of $Z$; see, e.g.,~\cite{beasley}. 

The RCB is linked to the nonnegative rank via an exact nonnegative factorization of $X$.
Indeed, suppose that $X=\sum_{k=1}^r W(:,k)H(k,:)$ is an exact NMF of $X$ with $r=\rank_+(X)$.
Let $\widehat{W}=\bin(W)$ and $\widehat{H}=\bin(H)$.
Then
\[
Z=\bin(X) \leq \text{bin}\left(\sum_{k=1}^rW(:,k)H(k,:)\right)\leq \sum_{k=1}^r\widehat{W}(:,k)\widehat{H}(k,:)=\widehat{W}\widehat{H}.
\]
Since $\widehat{W}\widehat{H}$ is a sum of $r$ binary rank-one rectangles that cover $\supp(Z)$, we have $\mathrm{RCB}(X)\leq \rank_+(X)$.

Note that $\mathrm{RCB}(X)\geq \mathrm{FSB}(X)$ for any nonnegative matrix $X$.
Indeed, two entries belonging to the same fooling set cannot lie in the same rectangle, so every rectangle cover must use at least one rectangle per entry in the fooling set.

\subsection{Standard approach via enumeration}

The classical approach to compute the RCB is to enumerate all maximal rectangles in $\mathcal{R}_{01}(Z)$ and solve a set-cover problem; see, e.g., \cite{fiorini2019, gillis2020nonnegative}. A maximal rectangle is a rectangle that is not contained in a larger rectangle. 
Let $R^{(1)},\dots,R^{(L)}$ be the maximal rectangles, and let the binary variable $z_\ell\in\{0,1\}$ indicate whether rectangle $\ell$ is selected.
This leads to the following integer program:
\begin{equation}\label{eq:rcb-enumFio}
\min_{z\in\{0,1\}^L}\quad \sum_{\ell=1}^L z_\ell 
\quad  \text{ such that } \quad 
 \sum_{\ell=1}^L z_\ell\,R^{(\ell)}\ \geq\ Z.
\end{equation}
This model uses $L$ binary variables and $|\supp(X)|$ covering constraints (the inequalities only need to be checked on the support of $X$). 
In general, the number of maximal rectangles, $L$, is exponential in the dimension. 
One can, for example, generate them by selecting all subsets of rows and then selecting the columns only when the corresponding entries in $X$ are positive. This leads to at most $L\leq 2^m$ rectangles. Using the same argument for the columns, we have $L\leq 2^{\min(m,n)}$, although there exist dedicated algorithms for this task, e.g.,~\cite{alexe2004consensus}. 

\subsection{Alternative approach using MILP}\label{subsec:RCBsecondapproach}

Another approach, introduced by Dewez~\cite{dewezphd}, checks directly whether $Z$
can be covered by at most $r$ rectangles without enumerating explicitly. 
It introduces $r$ binary matrices, 
$R^{(k)}\in\{0,1\}^{m\times n}$ for $k \in [r]$, 
where $R^{(k)}$ represents the support of the $k$th rectangle. 
The formulation enforces three conditions: (i)~every positive entry of $Z$ must be
covered by at least one of the matrices $R^{(k)}$, (ii)~no entry outside the support of $Z$ can be covered, 
and 
(iii)~each matrix $R^{(k)}$ must have a rectangular support.
This can be formulated as follows: 
\begin{equation}\label{eq:rcb-dewez}
\begin{aligned}
\text{Find}\quad & R^{(k)}\in\{0,1\}^{m\times n} \text{ for } k \in [r]\\
\text{s.t.}\quad
& \sum_{k=1}^r R^{(k)}\ge Z, \; 
R^{(k)} \le Z \text{ for all } 
k \in [r],\\
& R^{(k)}_{i_1j_1}+R^{(k)}_{i_2j_2}\le 1+R^{(k)}_{i_1j_2} \text{ for all } i_1\neq i_2, j_1\neq j_2,  k \in [r],\\
& R^{(k)}_{i_1j_1}+R^{(k)}_{i_2j_2}\le 1+R^{(k)}_{i_2j_1} \text{ for all } i_1\neq i_2, j_1\neq j_2, k \in [r].
\end{aligned}
\end{equation}
The first constraint ensures coverage of the support of $Z$; 
the second forbids entries outside the support, and 
the last constraint enforces the rectangle property for each $R^{(k)}$. 
This formulation uses $mnr$ binary variables and $\mathcal{O}(m^2n^2r)$ linear constraints.

\subsection{Proposed non-convex formulation}

As in the approach described in Section~\ref{subsec:RCBsecondapproach}, we seek a cover of the support of $X$ by at most $r$ rectangles without enumerating them.
Let $Z=\bin(X)$ and introduce binary matrices $U\in\{0,1\}^{m\times r}$ and $V\in\{0,1\}^{r\times n}$.
The $k$th rectangle is given by $R^{(k)}=U(:,k)V(k,:)$.
The proposed feasibility formulation is
\begin{equation}\label{eq:rcb-ncvx}
\begin{aligned}
\text{Find}\quad & U\in\{0,1\}^{m\times r},\quad V\in\{0,1\}^{r\times n}\\
\text{such that}\quad & Z\leq UV\leq rZ,
\end{aligned}
\end{equation}
where the inequalities are entrywise.
The $(i,j)$th entry of $UV$ counts the number of rectangles covering position $(i,j)$, because $U(i,k)V(k,j)$ equals $1$ if the $k$th rectangle covers this position and $0$ otherwise.
Hence $(UV)_{ij}\in\{0,\ldots,r\}$.
The inequalities $Z\leq UV\leq rZ$ ensure that the rectangles cover exactly the support of $Z$:
\begin{itemize}
    \item If $Z_{ij}=1$, then $1\leq(UV)_{ij}\leq r$. The lower bound requires at least one covering rectangle; the upper bound holds automatically because there are only $r$ rectangles. Overlaps are therefore allowed.
    \item If $Z_{ij}=0$, then $0\leq(UV)_{ij}\leq0$. Thus, no rectangle may cover this position. Since $U$ and $V$ are nonnegative, this is equivalent to $U(i,k)V(k,j)=0$ for every $k$.
\end{itemize}
The model uses $mr+rn$ binary variables, $|\supp(X)|$ bilinear covering constraints, and $r\bigl(mn-|\supp(X)|\bigr)$ bilinear admissibility constraints.
In particular, when $X$ is dense, the number of admissibility constraints remains small, which makes the formulation especially attractive.
Any feasible solution $(U,V)$ defines a rectangle cover of size at most $r$, and conversely, any such cover can be encoded in this way.
Figure~\ref{fig:rcledm} illustrates a cover of $M_5$.
\begin{figure}[!htbp]
\centering
\begin{tikzpicture}[scale=0.9, every node/.style={font=\small}]
    \def\s{0.42} 

    \begin{scope}[shift={(0,0)}]
    \node[font=\small\bfseries] at (1.05,-2.5) {$Z=\bin(M_5)$};

    \foreach \i in {0,...,4}{
        \foreach \j in {0,...,4}{
            \pgfmathtruncatemacro{\r}{\i+1}
            \pgfmathtruncatemacro{\c}{\j+1}
            \draw[black!45] (\j*\s,-\i*\s) rectangle ++(\s,-\s);
            \ifnum\r=\c
                \node at (\j*\s+0.5*\s,-\i*\s-0.5*\s) {\scriptsize 0};
            \else
                \node at (\j*\s+0.5*\s,-\i*\s-0.5*\s) {\scriptsize 1};
            \fi
        }
    }

    \foreach \r in {1,2,3}{
        \foreach \c in {4,5}{
            \fill[orange!45, opacity=0.45] ({(\c-1)*\s},{-(\r-1)*\s}) rectangle ++(\s,-\s);
        }
    }

    \foreach \r in {1,4,5}{
        \foreach \c in {2,3}{
            \fill[green!50!black!35, opacity=0.40] ({(\c-1)*\s},{-(\r-1)*\s}) rectangle ++(\s,-\s);
        }
    }

    \foreach \r in {2,4}{
        \foreach \c in {1,3,5}{
            \fill[magenta!45, opacity=0.35] ({(\c-1)*\s},{-(\r-1)*\s}) rectangle ++(\s,-\s);
        }
    }

    \foreach \r in {3,5}{
        \foreach \c in {1,2,4}{
            \fill[cyan!45, opacity=0.40] ({(\c-1)*\s},{-(\r-1)*\s}) rectangle ++(\s,-\s);
        }
    }

    \foreach \i in {0,...,4}{
        \foreach \j in {0,...,4}{
            \pgfmathtruncatemacro{\r}{\i+1}
            \pgfmathtruncatemacro{\c}{\j+1}
            \draw[black!45] (\j*\s,-\i*\s) rectangle ++(\s,-\s);
            \ifnum\r=\c
                \node at (\j*\s+0.5*\s,-\i*\s-0.5*\s) {\scriptsize 0};
            \else
                \node at (\j*\s+0.5*\s,-\i*\s-0.5*\s) {\scriptsize 1};
            \fi
        }
    }

    \foreach \i in {1,...,5}{
        \node[left]  at (-0.08,{-((\i-1)*\s+\s/2)}) {\scriptsize \i};
        \node[above] at ({(\i-1)*\s+\s/2},0.08) {\scriptsize \i};
    }
\end{scope}

    \begin{scope}[shift={(4.0,0)}]
        \node[font=\small\bfseries] at (0.84,-2.5) {$U$};

        \foreach \i in {0,...,4}{
            \foreach \j in {0,...,3}{
                \draw[black!45] (\j*\s,-\i*\s) rectangle ++(\s,-\s);
            }
        }

        \foreach \r in {1,2,3}{
            \fill[orange!55] (0,{-(\r-1)*\s}) rectangle ++(\s,-\s);
            \node at (0.5*\s,{-(\r-1)*\s-0.5*\s}) {\scriptsize 1};
        }
        \foreach \r in {4,5}{
            \fill[white] (0,{-(\r-1)*\s}) rectangle ++(\s,-\s);
            \node at (0.5*\s,{-(\r-1)*\s-0.5*\s}) {\scriptsize 0};
        }

        \foreach \r in {1,4,5}{
            \fill[green!50!black!35] (\s,{-(\r-1)*\s}) rectangle ++(\s,-\s);
            \node at (1.5*\s,{-(\r-1)*\s-0.5*\s}) {\scriptsize 1};
        }
        \foreach \r in {2,3}{
            \fill[white] (\s,{-(\r-1)*\s}) rectangle ++(\s,-\s);
            \node at (1.5*\s,{-(\r-1)*\s-0.5*\s}) {\scriptsize 0};
        }

        \foreach \r in {2,4}{
            \fill[magenta!45] (2*\s,{-(\r-1)*\s}) rectangle ++(\s,-\s);
            \node at (2.5*\s,{-(\r-1)*\s-0.5*\s}) {\scriptsize 1};
        }
        \foreach \r in {1,3,5}{
            \fill[white] (2*\s,{-(\r-1)*\s}) rectangle ++(\s,-\s);
            \node at (2.5*\s,{-(\r-1)*\s-0.5*\s}) {\scriptsize 0};
        }

        \foreach \r in {3,5}{
            \fill[cyan!45] (3*\s,{-(\r-1)*\s}) rectangle ++(\s,-\s);
            \node at (3.5*\s,{-(\r-1)*\s-0.5*\s}) {\scriptsize 1};
        }
        \foreach \r in {1,2,4}{
            \fill[white] (3*\s,{-(\r-1)*\s}) rectangle ++(\s,-\s);
            \node at (3.5*\s,{-(\r-1)*\s-0.5*\s}) {\scriptsize 0};
        }

        \foreach \i in {1,...,5}{
            \node[left] at (-0.08,{-((\i-1)*\s+\s/2)}) {\scriptsize \i};
        }
        \foreach \j in {1,...,4}{
            \node[above] at ({(\j-1)*\s+\s/2},0.08) {\scriptsize \j};
        }
    \end{scope}

    \begin{scope}[shift={(6.7,0)}]
        \node[font=\small\bfseries] at (1.05,-2.5) {$V$};

        \foreach \i in {0,...,3}{
            \foreach \j in {0,...,4}{
                \draw[black!45] (\j*\s,-\i*\s) rectangle ++(\s,-\s);
            }
        }

        \foreach \c in {4,5}{
            \fill[orange!55] ({(\c-1)*\s},0) rectangle ++(\s,-\s);
            \node at ({(\c-1)*\s+0.5*\s},-0.5*\s) {\scriptsize 1};
        }
        \foreach \c in {1,2,3}{
            \fill[white] ({(\c-1)*\s},0) rectangle ++(\s,-\s);
            \node at ({(\c-1)*\s+0.5*\s},-0.5*\s) {\scriptsize 0};
        }

        \foreach \c in {2,3}{
            \fill[green!50!black!35] ({(\c-1)*\s},-\s) rectangle ++(\s,-\s);
            \node at ({(\c-1)*\s+0.5*\s},-1.5*\s) {\scriptsize 1};
        }
        \foreach \c in {1,4,5}{
            \fill[white] ({(\c-1)*\s},-\s) rectangle ++(\s,-\s);
            \node at ({(\c-1)*\s+0.5*\s},-1.5*\s) {\scriptsize 0};
        }

        \foreach \c in {1,3,5}{
            \fill[magenta!45] ({(\c-1)*\s},-2*\s) rectangle ++(\s,-\s);
            \node at ({(\c-1)*\s+0.5*\s},-2.5*\s) {\scriptsize 1};
        }
        \foreach \c in {2,4}{
            \fill[white] ({(\c-1)*\s},-2*\s) rectangle ++(\s,-\s);
            \node at ({(\c-1)*\s+0.5*\s},-2.5*\s) {\scriptsize 0};
        }

        \foreach \c in {1,2,4}{
            \fill[cyan!45] ({(\c-1)*\s},-3*\s) rectangle ++(\s,-\s);
            \node at ({(\c-1)*\s+0.5*\s},-3.5*\s) {\scriptsize 1};
        }
        \foreach \c in {3,5}{
            \fill[white] ({(\c-1)*\s},-3*\s) rectangle ++(\s,-\s);
            \node at ({(\c-1)*\s+0.5*\s},-3.5*\s) {\scriptsize 0};
        }

        \foreach \i in {1,...,4}{
            \node[left] at (-0.08,{-((\i-1)*\s+\s/2)}) {\scriptsize \i};
        }
        \foreach \j in {1,...,5}{
            \node[above] at ({(\j-1)*\s+\s/2},0.08) {\scriptsize \j};
        }
    \end{scope}
\end{tikzpicture}
\caption{Illustration of 
    the RCB computation using~\eqref{eq:rcb-ncvx}, with  $\mathrm{RCB}(M_5)=4$. Left: the support matrix $Z=\bin(M_5)$, whose positive entries can be covered by four overlapping rectangles. Right: binary matrices $U$ and $V$ encoding the row and column sets of these four rectangles.}
\label{fig:rcledm}
\end{figure}

\section{Hyperplane separation bound (HSB)}\label{sec:hsb}

The HSB is a non-combinatorial lower bound on the nonnegative rank that relies on the magnitude of the entries of the input matrix $X$, and not only on its zero--nonzero pattern, unlike the two previous bounds.
It was popularized through the work of Rothvoß~\cite{rothvoss2017matching},  where it was used to show that the extension complexity of the matching polytope is exponential. 

In~\cite{rothvoss2017matching}, Rothvoß presents an unpublished result by Fiorini:  given a linear functional $L\in\mathbb{R}^{m\times n}$, the optimization problem \begin{equation}\label{alphaL}
\max_{u,v} ~\langle L,uv^\top \rangle \quad \text{such that }u\in[0,1]^m, v\in[0,1]^n,
\end{equation}
admits an optimal solution $(u^{\ast},v^{\ast})$ where $u^{\ast}\in\{0,1\}^m$ and $v^{\ast}\in\{0,1\}^n$.
By denoting the optimal objective of this problem $\alpha(L):= \langle L,u^{\ast}{v^{\ast}}^\top \rangle$ for a given $L\in\mathbb{R}^{m\times n}$, this result allows us to derive a lower bound on the nonnegative rank.

Assume $X\neq0$ and write $\|X\|_\infty=\max_{i,j}X_{ij}$.
Let $X=\sum_{k=1}^r W(:,k)H(k,:)$ be an exact NMF with $r=\operatorname{rank}_+(X)$.
Since $0\leq W(:,k)H(k,:) \leq X$, the rank-one matrix $\frac{1}{\|X\|_{\infty}}W(:,k)H(k,:)$ belongs to $[0,1]^{m\times n}$ for all $k$, hence
\[
\langle L,X\rangle = \|X\|_{\infty}\sum_{k=1}^r\left\langle L,\frac{W(:,k)H(k,:)}{\|X\|_{\infty}}\right\rangle\leq \|X\|_{\infty} \sum_{k=1}^r \alpha(L) = r\|X\|_{\infty} \alpha(L).
\]
For $\alpha(L)>0$, the above inequality gives
\[
\operatorname{rank}_+(X)\geq
\frac{\langle L,X/\|X\|_\infty\rangle}{\alpha(L)}.
\]
Maximizing over all such $L$ gives the hyperplane separation bound:
\begin{equation}\label{definitionHSB}
\mathrm{HSB}(X):= \max_{\substack{L\in\mathbb{R}^{m\times n}\\ \alpha(L)>0}} \frac{\left\langle L,X/\|X\|_\infty\right\rangle}{\alpha(L)}.
\end{equation}
The normalization by $\|X\|_\infty$ makes the bound invariant under positive scaling of $X$.
Since $\alpha(tL)=t\alpha(L)$ for every $t>0$, the ratio in~\eqref{definitionHSB} is unchanged when $L$ is multiplied by a positive scalar.
W.l.o.g, we may therefore normalize $\alpha(L)$ to one in~\eqref{definitionHSB}.
It follows from this observation that $\operatorname{HSB}(X)$ is the optimal value of the linear program that maximizes $\langle L,X/\|X\|_\infty\rangle$ subject to $\langle L,uv^\top\rangle\leq1$ for every $u\in\{0,1\}^m$ and $v\in\{0,1\}^n$.
These constraints are equivalent to imposing $\alpha(L)\leq 1$.
Any optimal solution of this LP satisfies $\alpha(L)=1$: indeed, since the objective $\langle L,X/|X|_\infty\rangle$ is positive at optimality (taking $L$ to be zero except for one entry equal to one at a position where $X$ attains its maximum gives a feasible solution with objective value one), replacing $L$ by $L/\alpha(L)$ would preserve feasibility and strictly increase the objective if $\alpha(L)<1$.

The LP described above imposes a constraint for every binary rank-one matrix.
To compute $\operatorname{HSB}(X)$, it suffices to retain the constraints associated with rectangles contained in the support of $X$:
\begin{equation}\label{eq:LP_HSB}
    \operatorname{HSB}(X)=
    \max_{L\in\mathbb{R}^{m\times n}}
    \left\langle L,\frac{X}{\|X\|_\infty}\right\rangle
    \quad\text{such that}\quad
    \langle L,R\rangle\leq1
    \quad\text{for all }R\in\mathcal{R}_{01}(X),
\end{equation}
where $\mathcal{R}_{01}(X)$ denotes the set of matrices $uv^\top$ with $u\in\{0,1\}^m$, $v\in\{0,1\}^n$, and $u_iv_j=0$ whenever $X_{ij}=0$.
Indeed, any feasible solution of this LP can be made feasible for the LP with all binary $u$ and $v$ by assigning sufficiently negative values to $L_{ij}$ whenever $X_{ij}=0$, without changing the objective.
This leaves the constraints indexed by $\mathcal{R}_{01}(X)$ unchanged.
All other constraints involve at least one entry $L_{ij}$ with $X_{ij}=0$ and are therefore satisfied if these entries are chosen sufficiently negative.
Conversely, every feasible solution of the LP with all binary $u$ and $v$ satisfies the constraints indexed by $\mathcal{R}_{01}(X)$.
The two LPs therefore have the same optimal value.
In the LP restricted to $\mathcal{R}_{01}(X)$, entries $L_{ij}$ with $X_{ij}=0$ appear neither in the objective nor in the constraints and may be fixed to zero.

\subsection{Existing implementations: enumeration and cutting planes}

Only a few works have focused on developing numerical tools to compute $\operatorname{HSB}(X)$.
Given a matrix $L$, a code is provided in~\cite{gillis2020nonnegative} to compute $\alpha(L)$ by solving~\eqref{alphaL} via the explicit enumeration of all binary rank-one matrices $uv^\top$.
Since the number of such factors can be as large as $2^{\min(m,n)}$, this approach does not scale well.

To the best of our knowledge, the only general implementation aimed at identifying the optimal matrix $L$ in~\eqref{definitionHSB} is that of~\cite{fiorini2019}, a conference paper in which the implementation details are not discussed.
However, an analysis of the associated code reveals a cutting-plane iterative scheme: starting from an initial pool of rank-one binary matrices, the algorithm alternates between (i)~solving the resulting LP to identify the best matrix $L$ under the constraints $\langle L,uv^\top\rangle \leq 1$ for all rank-one binary matrices $uv^\top$ currently in the pool, and 
(ii)~given $L$, searching for a rank-one matrix $uv^\top$ that violates the inequality $\langle L,uv^\top\rangle \leq 1$; this matrix is then added to the pool for the next LP.
In that implementation, the separation problem is formulated as a MILP by introducing $mn$ binary variables to linearize the products $u_iv_j$.

\subsection{Proposed non-convex formulation}

We adopt the same overall cutting-plane scheme as in~\cite{fiorini2019}, alternating between a master problem over $L$ and a separation problem over binary rank-one matrices.
\begin{itemize}
    \item We formulate the separation problem as a non-convex optimization problem and let the solver handle the bilinear terms.
More precisely, given a matrix $L\in\mathbb{R}^{m\times n}$, the separation problem consists in finding a matrix in $\mathcal{R}_{01}(X)$ that maximizes $\langle L,R\rangle$ by solving
\begin{equation}\label{eq:hsb-sep}
\max_{u\in\{0,1\}^m,\ v\in\{0,1\}^n} \ \langle L,uv^\top\rangle
\qquad
\text{ such that }\quad
u_i v_j = 0 \text{ for all } (i,j)\notin \supp(X).
\end{equation}
When solving~\eqref{eq:hsb-sep}, if a matrix $R\in\mathcal{R}_{01}(X)$ with $\langle L,R\rangle>1$ is found, the corresponding constraint $\langle L,R\rangle\leq1$ is added to the master problem.
If the optimal value of~\eqref{eq:hsb-sep} is at most one, it means that all the constraints indexed by $\mathcal{R}_{01}(X)$ are satisfied.
When $L$ is optimal for the current master problem, its objective value then equals $\operatorname{HSB}(X)$.
\item Let the finite set $\tilde{\mathcal R}\subseteq \mathcal{R}_{01}(X)$ denote the current pool of admissible binary rank-one matrices.
Our master problem is the following linear program 
\begin{equation}\label{eq:hsb-master}
\begin{aligned}
\max_{L\in\mathbb{R}^{m\times n}}\quad & \left\langle L,\frac{X}{\|X\|_\infty}\right\rangle\\
\text{such that}\quad & \langle L,R\rangle\leq1, & & \forall R\in\tilde{\mathcal{R}},\\
& L_{ij}\leq1, & & \forall (i,j)\in\supp(X),\\
& L_{ij}=0, & & \forall (i,j)\notin\supp(X).
\end{aligned}
\end{equation}

The singleton rectangles give the constraints $L_{ij}\leq1$ on $\operatorname{supp}(X)$.
With these constraints, the master problem is feasible and has a finite optimal value even when the pool $\tilde{\mathcal{R}}$ is empty; in that case, $L=\operatorname{bin}(X)$ is optimal.
\end{itemize}

Algorithm~\ref{algo:hsb} alternates separation~\eqref{eq:hsb-sep} and reoptimization of~\eqref{eq:hsb-master} after adding each generated constraint.
At any iteration, let $L$ be an optimal solution of~\eqref{eq:hsb-master}, and let $\bar{\alpha}$ be a finite positive upper bound on the optimal value of~\eqref{eq:hsb-sep}.
Then

\begin{equation}\label{eq:HSB-ineq}
\frac{\left\langle L,\frac{X}{\|X\|_\infty}\right\rangle}{\bar{\alpha}} \leq \mathrm{HSB}(X) \leq \left\langle L,\frac{X}{\|X\|_\infty}\right\rangle.
\end{equation}
Indeed, since $\langle L,R\rangle\leq\bar\alpha$ for every $R\in\mathcal{R}_{01}(X)$, $L/\bar\alpha$ satisfies all the constraints of the LP~\eqref{eq:LP_HSB}.
This proves the lower bound.
The upper bound follows because~\eqref{eq:hsb-master} retains only a subset of the rectangle constraints of that LP.
We use an upper bound $\bar{\alpha}$ since a global optimization solver can provide one without solving the separation problem to optimality; in our implementation (see Section~\ref{sec:numerics}), Gurobi returns this bound through its \texttt{ObjBound} attribute.

In practice, solving the separation problem~\eqref{eq:hsb-sep} to global optimality at every iteration is often unnecessary.
Indeed, in order to generate a new cut, it is sufficient to find any admissible binary rank-one matrix $R^+$ such that $\langle L,R^+\rangle>1$.
The separation search is first stopped if a feasible rectangle $R^+$ satisfies $\langle L,R^+\rangle>1+\delta$, where $\delta>0$ is a prescribed threshold. If no such rectangle is found, the search is continued until either a feasible rectangle with $\langle L,R^+\rangle>1$ is obtained or an upper bound $\bar\alpha\leq1+\epsilon$ is established, where $\epsilon>0$ is the required relative precision. This procedure is summarized in Algorithm~\ref{algo:hsb}.

\begin{algorithm}[h]
\caption{Cutting-plane algorithm for $\mathrm{HSB}(X)$}
\label{algo:hsb}
\begin{algorithmic}[1]
\Require $X\in\mathbb{R}_+^{m\times n}$, precision parameter $\epsilon$ (default = $10^{-6}$), threshold $\delta$ (default = $10^{-4}$)
\Ensure Bounds $B,U$ satisfying
$B\leq\mathrm{HSB}(X)\leq U$,
and the final master solution $L$

\State $X \leftarrow X/\|X\|_\infty$
\State $\tilde{\mathcal R}\leftarrow \emptyset$
\State $L \leftarrow \bin(X)$
\State $B \leftarrow 0$
\While{true}
    \State $U \leftarrow \langle L,X\rangle$
    \State Solve~\eqref{eq:hsb-sep}, stopping early if a feasible $R^+$ satisfies $\langle L,R^+\rangle > 1+\delta$
    \State Obtain $R^+$, if available,
    and an upper bound $\bar\alpha$ on the separation optimum
    \If{no such matrix $R^+$ has been found}
        \State Continue solving~\eqref{eq:hsb-sep} until $\bar\alpha\leq1+\epsilon$ or a feasible $R^+$ with $\langle L,R^+\rangle>1$
        is found
        \State Update $R^+$ and $\bar\alpha$
        \EndIf

    \If{$0<\bar\alpha<\infty$}
        \State $B\leftarrow
        \max\left\{B,\frac{U}{\bar\alpha}\right\}$
    \EndIf

    \If{$0<\bar\alpha\leq1+\epsilon$}
        \State \Return $B,U,L$
        \Comment{Requested precision attained}
    \EndIf

        \State $\tilde{\mathcal R}\leftarrow
        \tilde{\mathcal R}\cup\{R^+\}$
        \State Compute an optimal solution $L$
        of~\eqref{eq:hsb-master} for the updated pool
\EndWhile
\end{algorithmic}
\end{algorithm}

Let $U$ denote the final master objective and let $B$ be the largest lower bound obtained from~\eqref{eq:HSB-ineq} over the iterations. Since the final update gives $B\geq U/\bar\alpha$ and the stopping criterion ensures $\bar{\alpha}\leq 1+\epsilon$, we have
\[
B\leq\operatorname{HSB}(X)\leq U\leq(1+\epsilon)B.
\]
Thus, the returned lower and upper bounds satisfy $U/B-1\leq\epsilon$.

Figure~\ref{fig:hsconv} illustrates this behavior on $M_5$.
The solid curve shows the value of the master problem, which decreases as new rectangle constraints are added to the set $\tilde{\mathcal R}$, while the dashed curve shows the corresponding current lower bound. 
In this example, the final pool contains $27$ generated rectangles, and both bounds approach the exact value $\operatorname{HSB}(M_5)=2$, established in Proposition~\ref{prop:hsb-ledm}.
Note that this lower bound is rather weak since $\rank(M_5) = 3$; see Section~\ref{sec:numerics} for more details.
\begin{figure}[H]
    \centering






\begin{tikzpicture}
\begin{axis}[
    width=11.5cm,
    height=6.2cm,
    xmin=0, xmax=27,
    ymin=0.5, ymax=6.5,
    xtick={0,5,10,15,20,25,27},
    ytick={1,2,3,4,5,6},
    xlabel={Iterations (number of rank-one matrices in $\mathcal{\tilde{R}}$)},
    ylabel={Value},
    ymajorgrids=true,
    xmajorgrids=false,
    grid style={dotted, gray},
    axis y line*=left,
    axis x line*=bottom,
    tick align=outside,
    legend style={
        draw=none,
        fill=none,
        font=\small,
        at={(0.97,0.97)},
        anchor=north east
    },
]

\addplot[
    thick,
    blue!70!black,
    mark=*,
    mark size=1.8pt,
] coordinates {
    (0,6.2500000)
    (1,5.9375000)
    (2,5.9375000)
    (3,5.6250000)
    (4,5.3125000)
    (5,5.0000000)
    (6,4.6875000)
    (7,4.6875000)
    (8,4.3750000)
    (9,4.3750000)
    (10,4.3750000)
    (11,4.0000000)
    (12,3.6250000)
    (13,3.3750000)
    (14,3.2500000)
    (15,3.0000000)
    (16,2.6250000)
    (17,2.5000000)
    (18,2.1250000)
    (19,2.1250000)
    (20,2.1250000)
    (21,2.0625000)
    (22,2.0625000)
    (23,2.0625000)
    (24,2.0625000)
    (25,2.0000000)
    (26,2.0000000)
    (27,2.0000000)
};
\addlegendentry{$\langle L,X/\|X\|_\infty\rangle$}

\addplot[
    thick,
    orange!85!black,
    dashed,
    mark=square*,
    mark size=1.6pt,
] coordinates {
    (0,1.0416667)
    (1,0.9895833)
    (2,0.9895833)
    (3,0.9375000)
    (4,0.8854167)
    (5,0.8333333)
    (6,0.7812500)
    (7,0.7812500)
    (8,1.0937500)
    (9,1.0937500)
    (10,1.0937500)
    (11,1.0000000)
    (12,1.2083333)
    (13,1.1250000)
    (14,1.0833333)
    (15,1.0000000)
    (16,0.8750000)
    (17,0.8333333)
    (18,0.7083333)
    (19,0.7083333)
    (20,1.0625000)
    (21,1.0312500)
    (22,1.0312500)
    (23,1.0312500)
    (24,0.6875000)
    (25,1.0000000)
    (26,1.0000000)
    (27,2.0000000)
};
\addlegendentry{$\langle L,X/\|X\|_\infty\rangle/ \bar{\alpha}$}

\end{axis}
\end{tikzpicture}
    \caption{Illustration of the inequalities~\eqref{eq:HSB-ineq} on $M_5$. The solid curve shows the value $\langle L,X/\|X\|_\infty\rangle$ of the master problem, while the dashed curve shows the corresponding lower bound $\langle L,X/\|X\|_\infty\rangle/\bar{\alpha}$. 
    As binary rank-one matrices are added to the set $\tilde{\mathcal R}$, both curves approach the exact value $\operatorname{HSB}(M_5)=2$.}
    \label{fig:hsconv}
\end{figure}

\section{Self-scaled bound (SSB)} \label{sec:ssb}

The SSB, introduced by Fawzi and Parrilo~\cite{fawzi2016self}, follows the same
separation principle as the HSB: it searches for a linear functional~$L$ that is large on~$X$
while remaining small on a well-chosen set of  admissible  nonnegative rank-one matrices. 
The key difference lies in the admissible set.
While the HSB normalizes each rank-one factor by $\|X\|_{\infty}$ so that the admissible matrices lie in $[0,1]^{m\times n}$, the SSB instead requires each rank-one factor $R$ to satisfy $0\leq R\leq X$. 
This distinction has important computational consequences. In the HSB, the optimal solution of the separation problem is always attained at binary vectors, see Section~\ref{sec:hsb}, so the non-convexity reduces to a combinatorial problem over $\{0,1\}^{m\times n}$.
In the SSB, the variables remain continuous, and the constraint $0\leq uv^\top \leq X$ defines a non-convex feasible region. 
This makes the separation problem significantly harder to solve, which likely explains why the SSB has received less computational attention, although, as we will see, it leads to much stronger lower bounds.
The SSB is defined as
\begin{equation}\label{eq:ssb-def}
\mathrm{SSB}(X) \quad := \quad \max_{L \in \mathbb{R}^{m \times n}}
\langle L, X \rangle
~ \text{ such that } ~
\langle L, R \rangle \leq 1 \text{ for all }\, R \in \mathcal{R}(X), 
\end{equation}
where 
\[ \mathcal{R}(X) = \{uv^\top : u \in \mathbb{R}_+^m,\, v \in \mathbb{R}_+^n,\ 0\leq uv^\top\leq X\} 
\] 
is the set of admissible nonnegative rank-one matrices. 
To see that $\rank_+(X)\geq \mathrm{SSB}(X)$, let $X = \sum_{k=1}^r R_k$ be an exact NMF with $r = \rank_+(X)$.
Each factor $R_k$ is a nonnegative rank-one matrix satisfying $0 \leq R_k \leq X$, hence $R_k\in\mathcal{R}(X)$.
Therefore, for any feasible $L$, 
\[
\langle L, X \rangle = \sum_{k=1}^r \underbrace{\langle L, R_k\rangle}_{\leq 1} \leq r = \rank_+(X).
\]

It was shown in~\cite{fawzi2016self} that $\mathrm{SSB}(X)\geq \mathrm{HSB}(X)$ for every nonzero nonnegative matrix $X$. Indeed, the SSB keeps the full entrywise information of $X$ through the constraint $0\leq uv^\top \leq X$, whereas the HSB only uses the weaker condition $0\leq uv^\top \leq \|X\|_\infty \mathbf{1}\mathbf{1}^\top$ after normalization.
As a result, the maximization problem~\eqref{eq:ssb-def} defining the SSB is less constrained and leads to a larger value.
Fawzi and Parrilo~\cite{fawzi2016self} also proposed a computable SDP relaxation of the SSB.
Our approach directly evaluates the SSB itself through a cutting-plane method.
To the best of our knowledge, this is the first implementation that directly computes the SSB, rather than a relaxation of it.

\subsection{Proposed non-convex formulation}

We adopt the same cutting-plane scheme as for the HSB, alternating between a non-convex separation step and a linear master problem.
The main difference lies in the separation step, which now involves continuous variables.
Given a current matrix~$L$, instead of maximizing $\langle L, uv^\top \rangle$ over binary
vectors as for the HSB, we must solve the non-convex problem
\begin{equation}\label{eq:ssb-sep}
\max_{u \in \mathbb{R}_+^m,\, v \in \mathbb{R}_+^n}\;
\langle L, uv^\top \rangle
\quad \text{ such that } \quad
uv^\top \in\mathcal{R}(X).
\end{equation}
If the optimal value of the separation problem is at most $1$, all the constrains defining the SSB are satisfied.
If $L$ is optimal for~\eqref{eq:ssb-def}, its objective value then equals $\mathrm{SSB}(X)$.
In practice, global solvers handle such bilinear constraints through McCormick envelopes, which require explicit bounds on the variables.
Due to the scaling ambiguity in the factorization $uv^\top$, we impose without loss of generality that 
\[
0\leq u\le 1 
\quad \text{ and } \quad 
0\le v_j \le \max_i X_{ij} \text{ for } j=1,\dots,n. 
\] 
In fact, given any non-zero $uv^\top \leq X$, 
define $u' = \frac{u}{\|u\|_\infty} \in [0,1]^m$ with $\max_i u'_i = 1$, implying $v'_j = \|u\|_\infty v_j \le \max_i X_{ij}$, while $u' v'^\top = uv^\top$. 
The master problem has the same structure as in the HSB case: given the current set $\tilde{\mathcal R}$ of generated rank-one matrices, we solve
\begin{equation}\label{eq:ssb-master}
\begin{aligned}
\max_{L\in\mathbb{R}^{m\times n}}\quad & \langle L,X\rangle\\
\text{such that}\quad & \langle L,R\rangle\leq 1, & & \forall R\in\tilde{\mathcal R},\\
    & L_{ij}\leq \frac{1}{X_{ij}}, & & \forall (i,j)\in\supp(X),\\
    & L_{ij}=0, & & \forall (i,j)\notin\supp(X).
\end{aligned}
\end{equation}

For each $(i,j)\in\supp(X)$, the singleton $X_{ij}e_i e_j^\top$ belongs to $\mathcal R(X)$, and therefore any feasible $L$ satisfies $X_{ij}L_{ij}\leq1$.
Hence, $L_{ij}\leq 1/X_{ij}$ is a valid constraint for the SSB formulation.
Entries of $L$ outside $\supp(X)$ can be fixed to zero because every admissible rank-one matrix vanishes there and these entries do not contribute to the objective.
We include these constraints explicitly in the master problem.
They ensure that~\eqref{eq:ssb-master} is feasible and has a finite optimal value: $L=0$ is feasible, while
\[
\langle L,X\rangle
=\sum_{(i,j)\in\supp(X)}L_{ij}X_{ij}
\leq \sum_{(i,j)\in\supp(X)}\frac{1}{X_{ij}}X_{ij} =|\supp(X)|.
\]
When the pool is empty, this upper bound is attained by setting $L_{ij}=1/X_{ij}$ for $(i,j)\in\supp(X)$ and $L_{ij}=0$ otherwise.
This solution is therefore optimal for the initial master problem and is used to initialize Algorithm~\ref{algo:ssb}.

At any iteration, let $L$ be an optimal solution of the current master problem, and let $\bar\alpha$ be a finite positive upper bound on the optimal value of the separation problem~\eqref{eq:ssb-sep}. Then
\begin{equation}\label{eq:SSB-ineq}
\frac{\langle L,X\rangle}{\bar\alpha}
\leq \operatorname{SSB}(X)
\leq \langle L,X\rangle.
\end{equation}
Indeed, $\langle L,R\rangle\leq\bar\alpha$ for every $R\in\mathcal R(X)$, so $L/\bar\alpha$ is feasible for~\eqref{eq:ssb-def}, which proves the lower bound. The upper bound follows because the master problem retains only a subset of the constraints of~\eqref{eq:ssb-def}.

Following the same strategy as for the HSB, we rely on the global upper bound $\bar\alpha$ returned by the solver, which is available even when the separation problem is not solved to optimality.
Accordingly, the separation is performed in two stages: the search is first stopped when a feasible $R^+$ satisfies $\langle L,R^+\rangle>1+\delta$; if no such matrix is found, it continues until either a violated matrix with $\langle L,R^+\rangle>1$ is obtained or a global upper bound $\bar\alpha\leq1+\epsilon$ is established.
This procedure is summarized in Algorithm~\ref{algo:ssb}.
\begin{algorithm}[H]
\caption{Cutting-plane algorithm for $\mathrm{SSB}(X)$}
\label{algo:ssb}
\begin{algorithmic}[1]
\Require $X\in\mathbb{R}_+^{m\times n}$,
precision parameter $\epsilon>0$ (default $10^{-6}$), threshold $\delta$ (default $10^{-2}$)
\Ensure Bounds $B,U$ satisfying
$B\leq\mathrm{SSB}(X)\leq U$,
and the final master solution $L$

\State $\tilde{\mathcal R}\leftarrow\emptyset$
\State initialize $L_{ij}=1/X_{ij}$ for $(i,j)\in\supp(X)$, and $L_{ij}=0$ otherwise
\State $B\leftarrow0$

\While{true}
\State $U\leftarrow\langle L,X\rangle$

\State Solve~\eqref{eq:ssb-sep}, stopping early if a feasible $R^+$ satisfies
$\langle L,R^+\rangle>1+\delta$
\State Obtain $R^+$, if available, and an upper bound $\bar\alpha$ on the separation optimum

\If{no such matrix $R^+$ has been found}
    \State Continue solving~\eqref{eq:ssb-sep} until either
    $\bar\alpha\leq1+\epsilon$ or a feasible $R^+$ with
    $\langle L,R^+\rangle>1$ is found
    \State Update $R^+$ and $\bar\alpha$
\EndIf

\If{$0<\bar\alpha<\infty$}
    \State $B\leftarrow
    \max\left\{B,\frac{U}{\bar\alpha}\right\}$
\EndIf

\If{$0<\bar\alpha\leq1+\epsilon$}
    \State \Return $B,U,L$
\EndIf

\State $\tilde{\mathcal R}\leftarrow
\tilde{\mathcal R}\cup\{R^+\}$
\State Compute an optimal solution $L$ of~\eqref{eq:ssb-master}
for the updated pool

\EndWhile
\end{algorithmic}
\end{algorithm}
Let $U$ denote the final master objective and let $B$ be the largest lower bound obtained from~\eqref{eq:SSB-ineq} over the iterations. At termination, the stopping criterion $\bar\alpha\leq1+\epsilon$ gives
\[
B\leq\operatorname{SSB}(X)\leq U\leq(1+\epsilon)B.
\]
Let us illustrate how this enclosure can be used.
We consider the matrix $M_{11}$ (see Section~\ref{sec:numerics-ssb}), for which the
numerical estimate of the SSB lies just above $6$.
This raises the question of whether the estimate is sufficiently accurate to justify a lower bound of $7$ on the nonnegative rank.
Applied to $M_{11}$, using Gurobi (see Section~\ref{sec:numerics}), Algorithm~\ref{algo:ssb}, returns the bounds
\[
B \leq \operatorname{SSB}(M_{11}) \leq U,
\]
where
\[
B \approx 6.003617654226720,\qquad U \approx 6.003618232823751, \qquad \bar\alpha \approx 1 + 9.6374725\times10^{-8}.
\]
The absolute and relative gaps satisfy
\[
U-B < 5.8\times10^{-7}, \qquad \frac{U-B}{B}<10^{-7}.
\]
Both bounds round to $6.004$ at three decimal places, so the same holds for $\operatorname{SSB}(M_{11})$.
More importantly, the lower bound itself exceeds $6$, which gives
\[
\operatorname{rank}_+(M_{11})\geq\lceil B\rceil=7.
\]
Thus, in this case, the numerical tolerance does not affect the resulting integer lower bound.

Figure~\ref{fig:ssb_ledm5} provides a complementary view of the behavior of Algorithm~\ref{algo:ssb} throughout the iterations.
For $M_5$, the master upper bound decreases as new rank-one matrices are added to the pool, while the best lower bound progressively improves.
The two bounds tighten around the value $4.186$.
In particular, the final lower bound satisfies $B>4$, which yields $\operatorname{rank}_+(M_5)\geq5$ and therefore recovers the known value $\operatorname{rank}_+(M_5)=5$.

\begin{figure}[!htbp]
        \centering
        \begin{tikzpicture}
\begin{axis}[
    width=11.5cm,
    height=6.2cm,
    xmin=0, xmax=29,
    ymin=0.5, ymax=20.5,
    xtick={0,5,10,15,20,25,29},
    ytick={1,4,10,15,20},
    xlabel={Iterations (number of rank-one matrices in $\mathcal{\tilde{R}}$)},
    ylabel={Value},
    ymajorgrids=true,
    xmajorgrids=false,
    grid style={dotted, gray},
    axis y line*=left,
    axis x line*=bottom,
    tick align=outside,
    legend style={
        draw=none,
        fill=none,
        font=\small,
        at={(0.97,0.97)},
        anchor=north east
    },
]

\addplot[
    thick,
    blue!70!black,
    mark=*,
    mark size=1.8pt,
] coordinates {
    (0,20.0000000)
    (1,15.6473785)
    (2,11.2947555)
    (3,11.2947546)
    (4,11.2947515)
    (5,10.8318759)
    (6,9.5331035)
    (7,8.5481356)
    (8,8.0075949)
    (9,7.5447188)
    (10,6.4989379)
    (11,5.8347771)
    (12,5.5009619)
    (13,5.1583583)
    (14,5.0786067)
    (15,4.7832512)
    (16,4.6933825)
    (17,4.4256725)
    (18,4.3352283)
    (19,4.2886949)
    (20,4.2745142)
    (21,4.2481198)
    (22,4.2230155)
    (23,4.2213600)
    (24,4.2040829)
    (25,4.1941182)
    (26,4.1919855)
    (27,4.1894815)
    (28,4.1894815)
    (29,4.1856952)
};
\addlegendentry{$\langle L,X\rangle$}

\addplot[
    thick,
    orange!85!black,
    dashed,
    mark=square*,
    mark size=1.6pt,
] coordinates {
    (0,3.7364848)
    (1,2.9233099)
    (2,2.1101345)
    (3,2.1101343)
    (4,2.4167077)
    (5,2.3073820)
    (6,2.0397725)
    (7,1.8209046)
    (8,1.7133635)
    (9,1.7956255)
    (10,1.3905601)
    (11,1.3335728)
    (12,1.3327819)
    (13,1.2239717)
    (14,1.3592739)
    (15,1.5457049)
    (16,1.3895539)
    (17,1.4249293)
    (18,1.8085652)
    (19,2.6714191)
    (20,3.6410910)
    (21,3.4497878)
    (22,3.5105720)
    (23,3.2570134)
    (24,3.3926986)
    (25,3.7093962)
    (26,4.0475232)
    (27,3.8227640)
    (28,4.1423720)
    (29,4.1856944)
};
\addlegendentry{$\langle L,X\rangle/ \bar{\alpha}$}

\end{axis}
\end{tikzpicture}
        \caption{Illustration of the inequalities~\eqref{eq:SSB-ineq} on $M_5$. The solid curve shows the value $\langle L,X\rangle$ of the master problem, while the dashed curve shows the corresponding lower bound $\langle L,X\rangle/\bar{\alpha}$. 
    As admissible rank-one matrices are added to the set $\tilde{\mathcal R}$, both curves approach the value 4.186.}
        \label{fig:ssb_ledm5}
    \end{figure}



\Needspace{7\baselineskip}
\section{Numerical results}\label{sec:numerics}

In this section, we evaluate the behavior of the various approaches to compute the four bounds discussed in this paper (namely, FSB, RCB, HSB and SSB).
To this end, we consider a benchmark covering several matrix families with different properties. These families are described in Section~\ref{sec:datasets}. 
The code used in this paper is available online from  
\begin{center}
\url{https://github.com/timothy-baeckelant/LowerBoundsNCvx}
\end{center} 
Our main computational baseline is \cite{fiorini2019}, whose code is available online from
\begin{center}
\url{https://bitbucket.org/matthias-walter/nonnegrank/src/master/}
\end{center}
The computational times reported in the tables of this section are measured in seconds. Unless stated otherwise, a time limit of one hour was imposed. In some cases, however, we allowed the algorithms to run longer. This was especially the case for the SSB since, to the best of our knowledge, the values of this bound were not known for the matrices under consideration. Therefore, we aimed to obtain as many new values as possible. When several methods were available for the same bound, namely for the HSB and the RCB, and one method was run beyond the one-hour limit on a given matrix, the other methods were also run with a comparable time limit on that matrix. For instance, since our non-convex method, Algorithm~\ref{algo:hsb}, certified the HSB of the slack matrix of the icosidodecahedron in 22333 seconds (see \Tableref{tab:hsb-all}), we also ran the MILP implementation of~\cite{fiorini2019} on this instance with a time limit of 7 hours, that is, 25200 seconds.
All computations were performed on a machine with an Intel$^\text{®}$ Core\texttrademark{} i5-1335U 1.30\,GHz CPU, 16\,GB RAM, running Ubuntu~22.04. We use Gurobi~11 as the optimization solver. Unless stated otherwise, we report wall-clock times as arithmetic means over 30 independent runs with distinct random seeds. Although the formulations are mathematically exact, the numerical computations are subject to the solver's tolerances.
Gurobi's presolve and heuristic routines induce slight seed-dependent variations; averaging yields a stable and fair summary of computational cost. For each method, we report: bound value, runtime (s), and—when relevant—the number of iterations (for instance, the number of rank-one witnesses for SSB).
In each table, bold entries denote (i)~the method with the shortest running time, and (ii)~any method whose running time is within 10\% of this minimum.

\subsection{Non-convex optimization solver: Gurobi} \label{sec:gurobi}

All the non-convex models considered in this paper are solved with Gurobi.
At a high level, when a bilinear term appears in the formulation, Gurobi introduces an auxiliary variable for this product and reformulates the model
so that the original nonlinear constraints become linear constraints coupled with bilinear relations. These bilinear relations are then relaxed
through McCormick envelopes~\cite{mccormick1976computability}, which provide a convex outer approximation over the current variable bounds. The resulting relaxations are embedded in
a spatial branch-and-bound framework, where the solver may branch not only on integer variables, as in MILP, but also on continuous variables
involved in the non-convex terms. 
After branching, the McCormick envelopes are recomputed locally on the smaller domains, which tightens the
relaxation. 

As a result, the formulations proposed in Sections~\ref{sec:fs}--\ref{sec:ssb} can be written directly, without introducing explicit linearization schemes.
For example,
the products in \eqref{eq:ssb-sep} involve continuous variables and can be written directly in
GurobiPy, as illustrated in Listing~\ref{lst:ssb-separation-code}.

\begin{lstlisting}[style=pythonstyle,
caption={Implementation of the continuous bilinear separation problem for the SSB using GurobiPy.},
label={lst:ssb-separation-code}]
# Variable bounds used to remove the scaling ambiguity
u = model.addMVar(m, lb=0.0, ub=1.0)
v = model.addMVar(n, lb=0.0, ub=np.max(X, axis=0))

# Admissibility constraints: u v^T <= X
for i in range(m):
    for j in range(n):
        model.addConstr(u[i] * v[j] <= X[i, j])

# Maximize <L, u v^T>
model.setObjective(
    gp.quicksum(L[i,j] * u[i] * v[j] for i in range(m) for j in range(n)), GRB.MAXIMIZE)

model.Params.NonConvex = 2
\end{lstlisting}

This implementation point of view is also relevant when interpreting the numerical results.
Since Gurobi's treatment of bilinear and nonlinear non-convex constraints continues to improve over time\footnote{For nonconvex MIQCPs, Gurobi reports a 54.7\% speed-up in version 13 over version 12 on its overall benchmark set (\(>1\) s), and a \(2.68\times\) speed-up on the hard instances (\(>100\) s); see the Gurobi 13.0 presentation, \emph{What's New in Gurobi 13.0}.}, one may reasonably expect further speed-ups for the approaches developed in this paper.
By contrast, the standard MILP formulations used as baselines in the literature already rely on explicit linearizations.
Their computational bottleneck directly depends on the formulation size itself, which often involves very large numbers of variables and constraints.

\subsection{Benchmark nonnegative matrices} \label{sec:datasets}

We use the benchmark\footnote{See also \url{https://sites.google.com/site/exactnmf/data-set}.} of \cite{Vandaele_2015}, which contains families of matrices with different properties, 
including differences in their rank, sparsity pattern, and entry magnitudes.


\begin{enumerate}
\item \textbf{Linear Euclidean distance matrices (LEDMs).}
For $n\in\mathbb{N}$, the matrix $M_n\in\mathbb{R}_+^{n\times n}$ is defined by $M_n(i,j)=(i-j)^2$.
These matrices are dense, with zeros only on the diagonal, and satisfy $\rank(M_n)=3$ for all $n\geq 3$.
They were used in~\cite{beasley} to provide fixed rank matrices whose nonnegative rank increases with $n$; see Section~\ref{sec:intro}.  

\item \textbf{Slack matrices of regular $n$-gons.}
For $n\geq 3$, we denote by $S_n$ the slack matrix of the regular $n$-gons. 
These matrices also satisfy $\rank(S_n)=3$ and are a classical benchmark for lower bounds since their nonnegative rank is of particular interest~\cite{Fiorini_2012}.   
Since each vertex of an $n$-gon is incident to exactly two facets, and each facet contains exactly two vertices, the matrix $S_n$ has exactly two zero entries per row and per column.

\item \textbf{Unique-disjointness (UDISJ) matrices.}
For $n\in\mathbb{N}$, let $U_n$ be a partially specified
$2^n\times 2^n$ matrix indexed by pairs of binary vectors
$a,b\in\{0,1\}^n$, with
\begin{equation}\label{eq:udisj-definition}
U_n(a,b)=
\begin{cases}
1 & \text{if } a^\top b=0,\\
0 & \text{if } a^\top b=1,\\
? & \text{if } a^\top b>1.
\end{cases}
\end{equation}
A nonnegative completion of $U_n$ assigns a nonnegative value to each unspecified entry while preserving all specified entries.

We define $\mathrm{FSB}(U_n)$ as the maximum number of entries specified as $1$ such that, for every pair, at least one of the two cross entries is specified as $0$. Unspecified entries can neither be selected nor certify
the fooling-set condition. We define $\mathrm{RCB}(U_n)$ as the minimum number of rectangles covering all entries specified as $1$ while avoiding every entry specified as $0$. Unspecified entries may be covered but need not be covered.

Both quantities provide lower bounds on the nonnegative rank of every nonnegative completion of $U_n$. We do not evaluate the HSB or SSB on these partial patterns, since their definitions require a fully specified matrix. UDISJ patterns are classical hard instances in extension complexity~\cite{braun2016common}.

\item \textbf{Slack matrices of the dodecahedron, 24-cell, icosidodecahedron, and cuboctahedron.}
These are fixed polyhedral slack matrices coming from classical highly symmetric polytopes.
They complement the parametric families above by providing structured benchmark instances of moderate size arising from polyhedral geometry.

\item \textbf{Slack submatrices associated with the correlation polytope.}
For $n\in\mathbb{N}$, we denote by $C_n\in\mathbb{R}_+^{2^n\times 2^n}$ the matrix indexed by pairs of binary vectors $a,b\in\{0,1\}^n$ and defined by

\[
C_n(a,b)=(1-a^\top b)^2.
\]

The matrix $C_n$ is a fully specified nonnegative completion of $U_n$~\eqref{eq:udisj-definition}: it preserves the specified entries and assigns the strictly positive value $(1-a^\top b)^2$ to every previously unspecified position.

Every fooling set of the partial pattern $U_n$ is also a fooling set of $C_n$. Conversely, every rectangle cover of $C_n$ covers all entries specified as $1$ in $U_n$ and avoids all entries specified as $0$. Consequently,
\[
\mathrm{FSB}(U_n)\leq\mathrm{FSB}(C_n),
\qquad
\mathrm{RCB}(U_n)\leq\mathrm{RCB}(C_n).
\]

Moreover, these matrices are directly related to the correlation polytope: if $A,B\subseteq[n]$ denote the supports of $a$ and $b$, then $(|A\cap B|-1)^2$ is the slack at the vertex indexed by $B$ of the valid inequality indexed by $A$. Thus $C_n$ is a submatrix of a generalized slack matrix of the correlation polytope, whose rows need not define facets.
They are known to be difficult instances for nonnegative-rank lower bounds since it has been conjectured that their nonnegative rank matches the upper bound $2^n$. They were used in~\cite{fiorini2012linear} to provide lower bounds for the extension complexity of the traveling salesman polytope. 
\end{enumerate}
The instances used in the numerical experiments are summarized in Table~\ref{tab:datasets}.

\begin{table}[t]
\centering
\footnotesize
\setlength{\tabcolsep}{4pt}
\renewcommand{\arraystretch}{1.05}
\begin{tabular}{lccc}
\toprule
Instance & Dimensions & $\rank(X)$ & $\mathrm{nnz}(X)$ \\
\midrule
$M_n$ & $n\times n$ & $3$ &
$n^2-n$\\[1mm]

$S_n$ & $n\times n$ & $3$ &
$n^2-2n$ \\[1mm]

Dodecahedron & $20\times 12$ & $4$ &
$180 \;\; (75\%)$ \\

24-cell & $24\times 24$ & $5$ &
$432 \;\; (75\%)$ \\

Icosidodecahedron & $30\times 32$ & $4$ &
$840 \;\; (87.5\%)$ \\

Cuboctahedron & $12\times 14$ & $4$ &
$120 \;\; (71.43\%)$ \\[1mm]

$U_n$ & $2^n\times 2^n$ & $-$ &
$3^n$ \\[2mm]

$C_n$ & $2^n\times 2^n$ & $1+\dfrac{n(n+1)}2$ &
$4^n-n3^{n-1}$\\
\bottomrule
\end{tabular}
\caption{Summary of the benchmark matrix families used in the numerical experiments.
The quantity $\mathrm{nnz}(X)$ is the number of nonzero entries of $X$; for $U_n$, this only counts the specified nonzero entries.}
\label{tab:datasets}
\end{table}

\FloatBarrier
\subsection{Fooling Set Bound (FSB)}\label{sec:numerics-fs}

Table~\ref{tab:fsb-all} compares the standard enumerative IP formulation \eqref{eq:fsb-enumFio} of~\cite{fiorini2019} with our non-convex feasibility model~\eqref{eq:fsb-nc} on the benchmark instances described in Section~\ref{sec:datasets}.
For each matrix, we report the total runtime required to determine the value of the bound, together with the resulting FSB value. 
To compute the FSB, our non-convex approach needs to find a fooling set of size $r$ and certify that a set of size $r+1$ does not exist. 
Hence, runtimes are displayed as $a+b$, where $a$ denotes the cumulative time needed to certify the existence of a fooling set for each target size $r=1,\dots,\mathrm{FSB}(X)$, and $b$ denotes the time needed for the final infeasibility test certifying that no fooling set of size $\mathrm{FSB}(X)+1$ exists.

\begin{table*}[!htbp]
\centering
\setlength{\tabcolsep}{4pt}
\renewcommand{\arraystretch}{1.05}

\caption{Comparison of the MILP approach from~\cite{fiorini2019} and our non-convex formulation for the fooling set bound (FSB). 
Reported times are in seconds.} 
\label{tab:fsb-all}

\begin{subtable}[t]{0.48\textwidth}
\centering
\footnotesize
\caption{Linear Euclidean distance matrices}
\label{tab:fsb-ledm}
\begingroup
\fontsize{8}{9}\selectfont
\begin{tabular}{lccc}
\toprule
Matrix & \makecell{Enum.IP~\\\cite{fiorini2019}} & \makecell{Our NonCvx\\model~\eqref{eq:fsb-nc}} & FSB \\
\midrule
$M_{10}$  & \textbf{0.048} & 0.029 + 0.247 & 3 \\
$M_{20}$  & 3.96           & \textbf{0.028 + 1.45} & 3 \\
$M_{30}$  & 27.5           & \textbf{0.050 + 5.74} & 3 \\
$M_{40}$  & 644            & \textbf{0.088 + 11.91} & 3 \\
$M_{50}$  & $>1$h        & \textbf{0.142 + 26.12} & 3 \\
$M_{100}$ & $>1$h        & \textbf{0.778 + 185} & 3 \\
$M_{200}$ & $>1$h        & \textbf{4.71 + 1065} & 3 \\
$M_{240}$ & $>1$h        & \textbf{6.12 + 3590} & 3 \\
\bottomrule
\end{tabular}
\endgroup
\end{subtable}
\hfill
\begin{subtable}[t]{0.48\textwidth}
\centering
\footnotesize
\caption{Slack matrices of regular $n$-gons}
\label{tab:fsb-ngon}
\begingroup
\fontsize{8}{9}\selectfont
\begin{tabular}{lccc}
\toprule
Matrix & \makecell{Enum.IP~\\\cite{fiorini2019}} & \makecell{Our NonCvx\\model~\eqref{eq:fsb-nc}} & FSB \\
\midrule
$S_{5}$   & \textbf{0.018} & 0.016 + 0.02 & 5 \\
$S_{10}$  & \textbf{0.045} & 0.05 + 1.73 & 4 \\
$S_{20}$  & \textbf{5.58}  & 0.068 + 17.2 & 4 \\
$S_{30}$  & \textbf{39.7}  & 0.088 + 227 & 4 \\
$S_{35}$  & 490            & \textbf{0.119 + 299} & 4 \\
$S_{38}$  & $>1$h        & \textbf{0.131 + 562} & 4 \\
$S_{40}$  & $>1$h        & \textbf{0.209 + 1217} & 4 \\
$S_{43}$  & $>1$h        & \textbf{0.331 + 3295} & 4 \\
\bottomrule
\end{tabular}
\endgroup
\end{subtable}

\vspace{0.5em}

\begin{subtable}[t]{0.48\textwidth}
\centering
\footnotesize
\caption{Polyhedral slack matrices}
\label{tab:fsb-extra-slack}
\begingroup
\fontsize{8}{9}\selectfont
\begin{tabular}{lccc}
\toprule
Matrix & \makecell{Enum.IP~\\\cite{fiorini2019}} & \makecell{Our NonCvx\\model~\eqref{eq:fsb-nc}} & FSB \\
\midrule
Dodecahedron      & \textbf{0.707} & 1.78 + 114.5 & 6 \\
24-cell           & \textbf{18.34} & 6.23 + $>1$h & 8 \\
Icosidodecahedron & \textbf{274}   & 19.9 + 3197 & 6 \\
Cuboctahedron     & \textbf{0.213} & 0.836 + 3.64 & 6 \\
\bottomrule
\end{tabular}
\endgroup
\end{subtable}
\hfill
\begin{subtable}[t]{0.48\textwidth}
\centering
\footnotesize
\caption{Unique-disjointness matrices}
\label{tab:fsb-udisj}
\begingroup
\fontsize{8}{9}\selectfont
\begin{tabular}{lccc}
\toprule
Matrix & \makecell{Enum.IP~\\\cite{fiorini2019}} & \makecell{Our NonCvx\\model~\eqref{eq:fsb-nc}} & FSB \\
\midrule
$U_2$ & \textbf{0.003} & 0.030 + 0.040 & 3 \\
$U_3$ & \textbf{0.003} & 0.251 + 0.319 & 6 \\
$U_4$ & \textbf{0.007} & 0.178 + 0.229 & 9 \\
\bottomrule
\end{tabular}
\endgroup
\end{subtable}

\vspace{0.5em}

\begin{subtable}[t]{0.56\textwidth}
\centering
\footnotesize
\caption{Slack submatrices of the correlation polytope}
\label{tab:fsb-corr}
\begingroup
\fontsize{8}{9}\selectfont
\begin{tabular}{lccc}
\toprule
Matrix & \makecell{Enum.IP~\\\cite{fiorini2019}} & \makecell{Our NonCvx\\model~\eqref{eq:fsb-nc}} & FSB \\
\midrule
$C_2$ & \textbf{0.002} & 0.048 + 0.062 & 3 \\
$C_3$ & \textbf{0.004} & 0.415 + 0.511 & 7 \\
$C_4$ & \textbf{0.271} & 4.23 + 351 & 10 \\
\bottomrule
\end{tabular}
\endgroup
\end{subtable}

\end{table*}
\begin{figure}[!htbp]
\centering
\begin{tikzpicture}
\begin{semilogyaxis}[
    width=0.70\textwidth,
    xlabel={Dimension ($n$)}, 
    ylabel={Computation time (s)}, 
    xmin=5, 
    xmax=240, 
    ymin=0.1, 
    ymax=50000, 
    grid=both,
    major grid style={black!50}, 
    legend cell align={left}, 
    legend style={
        at={(0.98,0.02)},
        anchor=south east,
        nodes={scale=0.85}
    }
]

\addplot[
    blue,
    mark=*,
    mark size=1pt
]
coordinates {
    (5,0.007)
    (6,0.01)
    (7,0.012)
    (8,0.017)
    (9,0.03)
    (10,0.048)
    (11,0.106)
    (12,0.17)
    (13,0.328)
    (14,0.556)
    (15,0.881)
    (16,1.738)
    (17,1.987)
    (18,2.319)
    (19,2.662)
    (20,3.964)
    (21,5.183)
    (22,12.471)
    (23,11.663)
    (24,8.857)
    (25,8.859)
    (26,20.394)
    (27,13.314)
    (28,16.647)
    (29,21.767)
    (30,27.492)
    (31,33.416)
    (32,41.513)
    (33,50.591)
    (34,246.38)
    (35,73.252)
    (36,96.89)
    (37,407.919)
    (38,664.418)
    (39,170.361)
    (40,643.893)
    (41,1274)
    (42,1784)
    (43,1881)
    (44,2410)
    (45,2632)
};
\addlegendentry{Enum.IP~\NoHyper\cite{fiorini2019}\endNoHyper}

\addplot[
    red,
    mark=*,
    mark size=0.5pt
]
coordinates {
    (5,0.134)
    (6,0.189)
    (7,0.212)
    (8,0.239)
    (9,0.264)
    (10,0.276)
    (11,0.412)
    (12,0.588)
    (13,0.697)
    (14,0.878)
    (15,1.013)
    (16,1.091)
    (17,1.105)
    (18,1.256)
    (19,1.389)
    (20,1.478)
    (21,1.873)
    (22,2.372)
    (23,2.802)
    (24,3.25)
    (25,3.704)
    (26,4.114)
    (27,4.489)
    (28,4.983)
    (29,5.277)
    (30,5.822)
    (31,6.319)
    (32,7.157)
    (33,7.659)
    (34,8.245)
    (35,8.81)
    (36,9.619)
    (37,10.118)
    (38,10.786)
    (39,11.158)
    (40,12.054)
    (41,12.33)
    (42,15.75)
    (43,15.597)
    (44,16.242)
    (45,21.939)
    (46,22.845)
    (47,23.93)
    (48,24.15)
    (49,25.592)
    (50,26.441)
    (51,34.844)
    (52,43.587)
    (53,52.841)
    (54,61.886)
    (55,71.416)
    (56,80.01)
    (57,90.94)
    (58,96.625)
    (59,105.931)
    (60,114.596)
    (61,114.485)
    (62,110.307)
    (63,108.877)
    (64,105.552)
    (65,102.838)
    (66,100.29)
    (67,100.118)
    (68,95.708)
    (69,93.582)
    (70,91.865)
    (71,104.329)
    (72,111.83)
    (73,125.999)
    (74,132.794)
    (75,148.406)
    (76,156.091)
    (77,170.126)
    (78,178.344)
    (79,190.057)
    (80,195.245)
    (81,207.494)
    (82,216.377)
    (83,228.24)
    (84,236.85)
    (85,249.094)
    (86,260.426)
    (87,266.987)
    (88,284.298)
    (89,293.307)
    (90,300.087)
    (91,291.486)
    (92,274.96)
    (93,268.598)
    (94,260.466)
    (95,242.581)
    (96,234.197)
    (97,217.67)
    (98,211.019)
    (99,195.83)
    (100,183.42)
    (101,195.291)
    (102,199.551)
    (103,215.089)
    (104,216.754)
    (105,231.61)
    (106,236.617)
    (107,249.983)
    (108,261.234)
    (109,262.665)
    (110,275.006)
    (111,284.052)
    (112,292.693)
    (113,297.359)
    (114,315.134)
    (115,317.694)
    (116,331.741)
    (117,338.728)
    (118,342.091)
    (119,358.165)
    (120,361.063)
    (121,377.065)
    (122,381.479)
    (123,395.023)
    (124,400.986)
    (125,410.426)
    (126,415.622)
    (127,432.183)
    (128,435.774)
    (129,440.772)
    (130,452.877)
    (131,450.954)
    (132,464.917)
    (133,480.535)
    (134,482.231)
    (135,497.492)
    (136,502.558)
    (137,505.355)
    (138,517.463)
    (139,531.996)
    (140,541.311)
    (141,549.82)
    (142,560.443)
    (143,569.312)
    (144,573.132)
    (145,592.804)
    (146,589.249)
    (147,599.684)
    (148,619.629)
    (149,626.485)
    (150,632.864)
    (151,626.404)
    (152,656.252)
    (153,659.869)
    (154,676.333)
    (155,662.519)
    (156,690.804)
    (157,680.309)
    (158,701.687)
    (159,696.655)
    (160,726.106)
    (161,733.888)
    (162,735.844)
    (163,739.898)
    (164,738.544)
    (165,766.338)
    (166,767.744)
    (167,784.923)
    (168,798.383)
    (169,810.826)
    (170,815.981)
    (171,797.483)
    (172,817.604)
    (173,838.694)
    (174,828.856)
    (175,849.441)
    (176,842.279)
    (177,856.009)
    (178,858.389)
    (179,894.47)
    (180,883.063)
    (181,896.185)
    (182,926.195)
    (183,926.959)
    (184,910.897)
    (185,924.551)
    (186,950.556)
    (187,957.748)
    (188,953.535)
    (189,989.368)
    (190,985.79)
    (191,991.567)
    (192,1002.588)
    (193,1017.112)
    (194,1009.026)
    (195,1021.338)
    (196,1022.348)
    (197,1030.098)
    (198,1070.731)
    (199,1071.036)
    (200,1069.302)
    (201,1301.97)
    (202,1547.577)
    (203,1777.555)
    (204,2050.764)
    (205,2285.434)
    (206,2569.404)
    (207,2782.322)
    (208,3003.612)
    (209,3226.681)
    (210,3474.625)
    (211,3562.618)
    (212,3532.088)
    (213,3545.228)
    (214,3597.646)
    (215,3612.348)
    (216,3504.13)
    (217,3568.801)
    (218,3513.497)
    (219,3480.446)
    (220,3586.579)
    (221,3526.023)
    (222,3536.731)
    (223,3567.555)
    (224,3603.796)
    (225,3576.759)
    (226,3613.169)
    (227,3529.54)
    (228,3605.99)
    (229,3520.377)
    (230,3631.513)
    (231,3596.603)
    (232,3532.095)
    (233,3639.659)
    (234,3611.64)
    (235,3546.265)
    (236,3543.331)
    (237,3589.969)
    (238,3522.158)
    (239,3551.154)
    (240,3585.286)
};
\addlegendentry{Our NonCvx model~\NoHyper\eqref{eq:fsb-nc}\endNoHyper}

\addplot[
    brown!80!black,
    dashed,
    thick,
    mark=none
]
coordinates {
    (5,3600)
    (240,3600)
};
\addlegendentry{Time limit (1 h)}

\addplot[only marks, blue, mark=triangle, mark size=2pt, mark options={fill=white}, forget plot] coordinates {(46,3600)};
\node[anchor=south west, font=\scriptsize, blue] at (axis cs:46,3600) {$>1$ h};

\end{semilogyaxis}
\end{tikzpicture}
\caption{Semilogarithmic runtime for computing the FSB on LEDMs. The dashed line marks 3600 seconds. The open triangle at $n=46$ denotes an enumerative run exceeding one hour and is plotted at that limit.}\label{fig:fsb_ledm5}
\end{figure}

\paragraph{Discussion} 

Table~\ref{tab:fsb-all} highlights significant differences between the various families of matrices in the benchmark.
On LEDMs and regular $n$-gon slack matrices, the non-convex feasibility model becomes progressively more effective as the size increases and overtakes the standard enumerative IP formulation for the larger instances. 
This is consistent with the structure of the formulation~\eqref{eq:fsb-nc}: when the optimal fooling-set size $r$ remains small, the model only involves $O(mr+rn)$ binary variables and $O(r^2+r(mn-|\supp(X)|))$ bilinear constraints, whereas the standard formulation must handle a quadratic number of pairwise conflicts in the support. Figure~\ref{fig:fsb_ledm5} provides a visual illustration of the significant difference in runtimes for LEDMs. 

The runtime decomposition also shows that, for the non-convex approach, most of the computational effort typically comes from proving that no fooling set of size $\mathrm{FSB}(X)+1$ exists.
In contrast, the successive feasibility checks for the target sizes $r=1,\dots,\mathrm{FSB}(X)$ are often much cheaper.

On the other benchmark matrices, the standard enumerative IP remains clearly faster on the instances considered here. This is likely due to the fact that these matrices have sparse supports, for which the explicit formulation of the conflicts appears to remain more manageable. In such cases, the advantage of our non-convex approach in terms of variable count is less decisive, while the cost of solving a sequence of global bilinear feasibility problems remains significant.
The final infeasibility test for size $\mathrm{FSB}(X)+1$ is typically the main bottleneck.

Overall, the experiments suggest that the proposed non-convex approach is most attractive on dense and structured matrices with small optimal fooling sets, whereas the standard enumerative IP remains preferable on sparse matrices. 

\FloatBarrier
\subsection{Rectangle Covering Bound (RCB)}\label{sec:numerics-rc}

Table~\ref{tab:rcb-all} compares the standard enumerative formulation of~\cite{fiorini2019}, the MILP formulation of~\cite{dewezphd}, and our non-convex feasibility model~\eqref{eq:rcb-ncvx}. 
For each matrix, we report the total runtime required to determine the value of the bound, together with the resulting RCB value.
For the two formulations with fixed $r$, tests start at $r=\min(m,n)$ and proceed in decreasing order. Runtimes are displayed as $a+b$, where $a$ is the cumulative time spent testing the feasible target sizes $\min(m,n),\ldots,\mathrm{RCB}(X)$, and $b$ denotes the time needed to certify that no rectangle cover of size $\mathrm{RCB}(X)-1$ exists.

\begin{table*}[!htbp]
\centering
\setlength{\tabcolsep}{3pt}
\renewcommand{\arraystretch}{1.05}

\caption{Comparison of the standard and non-convex formulations for the rectangle covering bound (RCB). Reported times are in seconds. 
For our approach and that of \cite{dewezphd}, 
runtimes are displayed as $a+b$, where $a$ is the cumulative time spent to find a feasibe rectangle cover if size $\mathrm{RCB}(X)$, and $b$ is the time spent to certify that no rectancle cover of size $\mathrm{RCB}(X)-1$ exists.}
\label{tab:rcb-all}

\begin{subtable}[t]{0.48\textwidth}
\centering
\footnotesize
\caption{Linear Euclidean distance matrices}
\label{tab:rcb-ledm}
\begingroup
\fontsize{7.5}{8.5}\selectfont
\begin{tabular}{@{}lcccc@{}}
\toprule
Matrix & \makecell{Enum.IP\\\cite{fiorini2019}} & \makecell{MILP\\\cite{dewezphd}} & \makecell{Our NonCvx\\model \eqref{eq:rcb-ncvx}} & RCB \\
\midrule
$M_{5}$   & \textbf{0.004} & 0.02 + 0.01 & 0.03 + 0.03 & 4 \\
$M_{7}$   & \textbf{0.036} & 0.11 + 0.74 & 0.05 + 0.06 & 5 \\
$M_{8}$   & 0.107          & 0.21 + 1.8   & \textbf{0.04 + 0.01} & 5 \\
$M_{10}$  & 0.959          & 0.71 + 11.8 & \textbf{0.10 + 0.58} & 5 \\
$M_{11}$  & 4.32           & 1.83 + 248.4  & \textbf{0.12 + 2.51} & 6 \\
$M_{15}$  & 985.8          & 11.9 + $>$1h & \textbf{0.24 + 3.92} & 6 \\
$M_{20}$  & $>$1h        & $>$1h       & \textbf{0.82 + 21.5} & 6 \\
$M_{21-35}$ & $>$1h & $>$1h & $<1$s + $>$1h & 7 \\
\bottomrule
\end{tabular}
\endgroup
\end{subtable}
\hfill
\begin{subtable}[t]{0.48\textwidth}
\centering
\footnotesize
\caption{Slack matrices of regular $n$-gons}
\label{tab:rcb-ngon}
\begingroup
\fontsize{7.5}{8.5}\selectfont
\begin{tabular}{@{}lcccc@{}}
\toprule
Matrix & \makecell{Enum.IP\\\cite{fiorini2019}} & \makecell{MILP\\\cite{dewezphd}} & \makecell{Our NonCvx\\model \eqref{eq:rcb-ncvx}} & RCB \\
\midrule
$S_{5}$  & \textbf{0.003} & 0.03 + 0.07  & 0.01 + 0.04 & 5 \\
$S_{7}$  & \textbf{0.009} & 0.15 + 0.77 & 0.02 + 0.29 & 6 \\
$S_{10}$ & \textbf{0.692} & 4.89 + 1836   & 0.09 + 1.78 & 7 \\
$S_{12}$ & \textbf{4.208} & 26.09 + $>$1h & 0.09 + 16.89 & 7 \\
$S_{13}$ & 17.24          & 253.4 + $>$1h & \textbf{0.11 + 14.13} & 7 \\
$S_{16}$ & $>$1h        & $>$1h       & \textbf{2.01 + 66.12} & 7 \\
$S_{19}$ & $>$1h        & $>$1h       & \textbf{18.1 + 122.9} & 7 \\
$S_{20}$ & $>$1h        & $>$1h       & $>$1h + 289.5 & 7 \\
\bottomrule
\end{tabular}
\endgroup
\end{subtable}

\vspace{0.5em}

\begin{subtable}[t]{0.48\textwidth}
\centering
\footnotesize
\caption{Polyhedral slack matrices}
\label{tab:rcb-extra-slack}
\begingroup
\fontsize{7.5}{8.5}\selectfont
\begin{tabular}{@{}lcccc@{}}
\toprule
Matrix & \makecell{Enum.IP\\\cite{fiorini2019}} & \makecell{MILP\\\cite{dewezphd}} & \makecell{Our NonCvx\\model \eqref{eq:rcb-ncvx}} & RCB \\
\midrule
Dodeca.      & \textbf{1.02} & $>$1h + 894 & 0.721 + 50.9 & 9 \\
24-cell           & \textbf{1516}  & $>$1h       & $>$1h & 12 \\
Icosidodec. & $>$1h            & $>$1h       & \revise{1672 + $>$1h} & \revise{$\geq 8$} \\
Cubocta.     & \textbf{0.02} & 553 + $>$1h & 0.23 + 1.11 & 8 \\
\bottomrule
\end{tabular}
\endgroup
\end{subtable}
\hfill
\begin{subtable}[t]{0.48\textwidth}
\centering
\footnotesize
\caption{Unique-disjointness matrices}
\label{tab:rcb-udisj}
\begingroup
\fontsize{7.5}{8.5}\selectfont
\begin{tabular}{@{}lcccc@{}}
\toprule
Matrix & \makecell{Enum.IP\\\cite{fiorini2019}} & \makecell{MILP\\\cite{dewezphd}} & \makecell{Our NonCvx\\model \eqref{eq:rcb-ncvx}} & RCB \\
\midrule
$U_2$ & \textbf{0.001} & 0.003 + 0.005 & 0.004 + 0.004 & 3 \\
$U_3$ & \textbf{0.003} & 0.173 + 0.157 & 0.007 + 0.053 & 6 \\
$U_4$ & \textbf{0.002} & 25 + 9.81     & 0.007 + 0.054 & 9 \\
$U_5$ & \textbf{1.352} & $>$1h     & $>$1h + $>$1h & 18 \\
$U_6$ & \textbf{1216} & $>$1h   & $>$1h & 27 \\

\bottomrule
\end{tabular}
\endgroup
\end{subtable}

\vspace{0.5em}

\begin{subtable}[t]{0.56\textwidth}
\centering
\footnotesize
\caption{Slack submatrices of the correlation polytope}
\label{tab:rcb-corr}
\begingroup
\fontsize{7.5}{8.5}\selectfont
\begin{tabular}{@{}lcccc@{}}
\toprule
Matrix & \makecell{Enum.IP\\\cite{fiorini2019}} & \makecell{MILP\\\cite{dewezphd}} & \makecell{Our NonCvx\\model \eqref{eq:rcb-ncvx}} & RCB \\
\midrule
$C_2$ & \textbf{0.002} & 0.003 + 0.005 & 0.003 + 0.005 & 3 \\
$C_3$ & \textbf{0.013} & 0.301 + 0.257 & 0.006 + 0.056 & 7 \\
$C_4$ & \textbf{0.162} & $>1$h       & 1.99 + 85.9 & 13 \\
$C_5$ & \textbf{2600} & $>1$h       &  $>1$h & 25\\
\bottomrule
\end{tabular}
\endgroup
\end{subtable}

\end{table*}

\paragraph{Discussion} 
The results show that the proposed non-convex model is particularly effective on the LEDMs and the slack matrices of the regular $n$-gons.
On both families, the standard enumerative approach is the fastest on small instances, but the non-convex formulation becomes superior as the dimension increases.
The transition occurs around $S_{13}$ for regular $n$-gons and already around $M_8$ for LEDMs, after which the gap becomes substantial.
The MILP formulation from~\cite{dewezphd} does not appear competitive in these experiments.

The runtime decomposition further shows that, for our non-convex approach, the dominant cost usually comes from the final feasibility test certifying that no rectangle cover of size $\mathrm{RCB}(X)-1$ exists. This is especially visible on LEDMs and regular $n$-gons, where the first term remains small while the second one accounts for most of the total runtime.

As for the FSB, the standard enumerative approach remains faster on the other benchmark matrices.
A plausible explanation is that for these sparser matrices, the number of feasible rectangles is smaller, whereas our non-convex approach requires solving a sequence of global bilinear feasibility problems.

\FloatBarrier
\subsection{Hyperplane separation bound (HSB)}\label{sec:numerics-hsb}

Table~\ref{tab:hsb-all} compares the MILP implementation of~\cite{fiorini2019}, along with a variant obtained by modifying the symmetry and enumeration parameters in the same code, and our non-convex approach on the benchmark instances described in Section~\ref{sec:datasets}.
For each matrix, we report the total runtime required to compute the bound, together with the resulting HSB value.

\begin{table*}[!htbp]
\centering
\captionsetup[subtable]{skip=2pt}
\setlength{\tabcolsep}{1.5pt}
\renewcommand{\arraystretch}{1.05}

\caption{Running time in seconds of the different methods to compute $\mathrm{HSB}(X)$.
Comparison of the standard and non-convex formulations for the $\mathrm{HSB}$. 
Comparison of the MILP implementation of~\cite{fiorini2019}, its variant with symmetry handling and enumeration, and our non-convex cutting-plane approach (Algorithm \ref{algo:hsb})}
\label{tab:hsb-all}

\begin{subtable}[t]{0.48\textwidth}
\centering
\footnotesize
\caption{Slack matrices of regular $n$-gons}
\label{tab:hsb-ngon}
\begingroup
\fontsize{7.5}{8.5}\selectfont
\begin{tabular}{@{}lcccc@{}}
\toprule
Matrix & \makecell{MILP\cite{fiorini2019}\\} & \makecell{MILP\cite{fiorini2019}\\+sym,enum} & \makecell{Our NonCvx\\meth.(Algo~\ref{algo:hsb})} & HSB \\
\midrule
$S_{5}$  & \textbf{0.014} & \textbf{0.014} & 0.106 & 3.090 \\
$S_6$ & \textbf{0.025} & \textbf{0.023} & 0.169 & 3.000\\
$S_7$ & \textbf{0.036} & \textbf{0.039} & 0.427 & 3.115 \\
$S_{8}$  & \textbf{0.066} & 0.089 & 1.742 & 3.314 \\
$S_9$ & \textbf{0.303} & \textbf{0.311} & 2.11 & 3.126 \\
$S_{10}$ & \textbf{0.48} & 0.69 & 5.58 & 3.090 \\
$S_{11}$ & \textbf{1.167} & 4.88 & 11.02 & 3.046\\
$S_{12}$ & \textbf{18.92} & 21.12 & 20.83 & 3.054 \\
$S_{13}$ & \textbf{11.64} & 21.02 & 34.37 & 3.043 \\
$S_{14}$ & \textbf{57.5} & 69.1 & 89.17 & 3.056 \\
$S_{15}$ & 94.4 & \textbf{86.23} & 117 & 3.000 \\
$S_{16}$ & 431 & 420 & \textbf{227} & 3.052 \\
$S_{17}$ & 381 & 361 & \textbf{323} & 3.018 \\
$S_{18}$ & 1011 & 1087 & \textbf{537} & 3.023 \\
$S_{19}$ & 2799 & 3071 & \textbf{515} & 3.025 \\
$S_{20}$ & $>1$h & $>1$h & \textbf{1214} & 3.020 \\
$S_{21}$ & $>1$h & $>1$h & \textbf{597} & 3.000 \\
$S_{22}$ & $>1$h & $>1$h & \textbf{1227} & 3.025 \\
$S_{23}$ & $>1$h & $>1$h & \textbf{3524} & 3.008 \\
$S_{24}$ & $>2$h & $>2$h & \textbf{5587} & 3.010 \\
\bottomrule
\end{tabular}
\endgroup
\end{subtable}
\hfill
\begin{subtable}[t]{0.48\textwidth}
\centering
\footnotesize
\caption{Linear Euclidean distance matrices}
\label{tab:hsb-ledm}
\begingroup
\fontsize{7.5}{8.5}\selectfont
\begin{tabular}{@{}lcccc@{}}
\toprule
Matrix & \makecell{MILP\cite{fiorini2019}\\} & \makecell{MILP\cite{fiorini2019}\\+sym,enum} & \makecell{Our NonCvx\\meth.(Algo~\ref{algo:hsb})} & HSB \\
\midrule
$M_5$ & \textbf{0.014} & \textbf{0.014} & 0.122 & 2\\
$M_6$ & \textbf{0.016} & 0.018 & 0.246 & 2\\
$M_7$ & \textbf{0.029} & 0.035 & 0.398 & 2 \\
$M_8$ & 0.045 & \textbf{0.037} & 0.789 & 2 \\
$M_9$ & 0.112 & \textbf{0.078} & 1.347 & 2 \\
$M_{10}$ & \textbf{0.213} & \textbf{0.198} & 2.754 & 2 \\
$M_{11}$ & \textbf{1.01} & \textbf{0.931} & 5.281 & 2\\
$M_{12}$ & \textbf{1.57} & 2.367 & 6.381 & 2 \\
$M_{13}$ & \textbf{2.92} & 8.779 & 8.646 & 2 \\
$M_{14}$ & \textbf{4.54} & 32.934 & 14.2 & 2 \\
$M_{15}$ & \textbf{8.28} & 126 & 13.2 & 2 \\
$M_{16}$ & 18.8 & 433 & \textbf{15.3} & 2 \\
$M_{17}$ & 32.3 & 1555 & \textbf{20.7} & 2 \\
$M_{18}$ & 88 & $>$1h & \textbf{52.8} & 2 \\
$M_{19}$ & 192 & $>$1h & \textbf{71.5} & 2\\
$M_{20}$ & 342 & $>$1h & \textbf{122} & 2\\
$M_{21}$ & 588 & $>$1h & \textbf{169} & 2\\
$M_{22}$ & 1346 & $>$1h & \textbf{87.8} & 2\\
$M_{23}$ & $>$1h & $>$1h & \textbf{149} & 2\\
$M_{24}$ & $>$1h & $>$1h & \textbf{174} & 2\\
$M_{25}$ & $>$1h & $>$1h & \textbf{240} & 2\\
$M_{26}$ & $>$1h & $>$1h & \textbf{361} & 2\\
$M_{27}$ & $>$1h & $>$1h & \textbf{1432} & 2\\
$M_{28}$ & $>$1h & $>$1h & \textbf{1964} & 2\\
$M_{29}$ & $>$1h & $>$1h & \textbf{2778} & 2\\
$M_{30}$ & $>$1h & $>$1h & \textbf{1946} & 2\\
\bottomrule
\end{tabular}
\endgroup
\end{subtable}

\vspace{0.5em}

\begin{subtable}[t]{0.48\textwidth}
\centering
\footnotesize
\caption{Polyhedral slack matrices}
\label{tab:hsb-extra-slack}
\begingroup
\fontsize{7.5}{8.5}\selectfont
\begin{tabular}{@{}lcccc@{}}
\toprule
Matrix & \makecell{MILP\cite{fiorini2019}\\} & \makecell{MILP\cite{fiorini2019}\\+sym,enum} & \makecell{Our NonCvx\\meth.(Algo~\ref{algo:hsb})} & HSB \\
\midrule
Dodeca.      & 54 & \textbf{6.346} & 332 & 3.819 \\
24-cell           & $>$1h & \textbf{5.091} & $>$1h & 5.600 \\
Icosidodeca. & $>$7h & -- & \textbf{22333} & 3.708 \\
Cubocta.     & 11.79 & \textbf{2.046} & 75.57 & 4.000 \\
\bottomrule
\end{tabular}
\endgroup
\end{subtable}
\hfill
\begin{subtable}[t]{0.48\textwidth}
\centering
\footnotesize
\caption{Slack submatrices of the correlation polytope}
\label{tab:hsb-corr}
\begingroup
\fontsize{7.5}{8.5}\selectfont
\begin{tabular}{@{}lcccc@{}}
\toprule
Matrix & \makecell{MILP\cite{fiorini2019}\\} & \makecell{MILP\cite{fiorini2019}\\+sym,enum} & \makecell{Our NonCvx\\meth.(Algo~\ref{algo:hsb})} & HSB \\
\midrule
$C_{2}$ & \textbf{0.015} & \textbf{0.014} & 0.164 & 3.667 \\
$C_{3}$ & \textbf{0.027} & 0.031 & 0.98 & 2.500 \\
    $C_{4}$ & 0.887 & \textbf{0.768} & 52.26 & 2.032 \\
$C_{5}$ & \textbf{7934} & $>$3h & $>$3h & 2.100 \\
\bottomrule
\end{tabular}
\endgroup
\end{subtable}

\end{table*}

\paragraph{A constant HSB on LEDMs}
The computations reported in Table~\ref{tab:hsb-all} show that $\mathrm{HSB}(M_n)=2$ for every LEDM instance considered. This led us to suspect that the same value might hold for the entire family and suggested the structure underlying the proof.
The following proposition confirms this observation for all $n$.

\begin{proposition}[HSB of linear Euclidean distance matrices]
\label{prop:hsb-ledm}
For every integer $n\geq2$,
\[
\operatorname{HSB}(M_n)=2.
\]
\end{proposition}

\begin{proof}
We use the linear programming formulation~\eqref{eq:LP_HSB}. For the lower bound, let $L_{1n}=L_{n1}=1$ and set all other entries of $L$ to zero. No rectangle in $\mathcal R_{01}(M_n)$ can contain both $(1,n)$ and $(n,1)$, since $M_{1,1}=M_{n,n}=0$. Hence,
\[
\langle L,R\rangle\leq1
\qquad\text{for every }R\in\mathcal R_{01}(M_n),
\]
so $L$ is feasible for~\eqref{eq:LP_HSB}. Since
\[
\langle L,M_n\rangle=2(n-1)^2
\qquad\text{and}\qquad
\|M_n\|_\infty=(n-1)^2,
\]
we obtain $\operatorname{HSB}(M_n)\geq2$.

For the reverse inequality, we decompose $M_n$ into binary rank-one matrices belonging to $\mathcal{R}_{01}(M_n)$.
For $1\leq k,\ell\leq n-1$, we define the matrix
$R^{k\ell}\in\{0,1\}^{n\times n}$ by
\[
R^{k\ell}_{ij}
=
\mathbf{1}_{\{i>\max(k,\ell)\}}
\mathbf{1}_{\{j\leq\min(k,\ell)\}},
\qquad 1\leq i,j\leq n.
\]
Thus, $R^{k\ell}$ has a rectangular block of ones in rows
$\max(k,\ell)+1,\ldots,n$ and columns $1,\ldots,\min(k,\ell)$, and zeros elsewhere; in particular, its support lies entirely in the strict lower triangle.
It means that each matrix $R^{k\ell}$ and its transpose belong to $\mathcal R_{01}(M_n)$. Moreover, for $i>j$, exactly $(i-j)^2$ pairs $(k,\ell)$ satisfy $j\leq k,\ell<i$. Therefore,
\[
M_n
=
\sum_{k,\ell=1}^{n-1}
\left(R^{k\ell}+(R^{k\ell})^\top\right).
\]
For any feasible $L$ in~\eqref{eq:LP_HSB}, each term in this sum has scalar product with $L$ at most one. Consequently,
\[
\langle L,M_n\rangle
\leq 2(n-1)^2
=
2\|M_n\|_\infty.
\]
It follows that $\operatorname{HSB}(M_n)\leq2$, which completes the proof.
\end{proof}

For $n\geq3$, this is weaker than the ordinary-rank bound $\rank(M_n)=3$; moreover, HSB remains constant while $\rank_+(M_n)$ grows unbounded. 

\paragraph{Discussion}
For regular $n$-gon slack matrices, the HSB stays close to $3$ although the values are not monotone in the dimension $n$.
Together with Proposition~\ref{prop:hsb-ledm} for the LEDMs, this suggests that HSB captures only a limited part of the complexity of these families.

From a computational point of view, the proposed non-convex approach becomes competitive on larger regular $n$-gons.
For example, on $S_{19}$, it requires $515$ seconds, compared with $2799$ seconds for the fastest competing implementation, corresponding to a speed-up factor of approximately $5.4$.
In general, the computational effort is large compared to the quality of the resulting bounds.
The numerical results are somewhat disappointing for the matrices considered in this benchmark, even though HSB is the bound that played a central role in Rothvo{\ss}' seminal work~\cite{rothvoss2017matching}. 

\FloatBarrier
\subsection{Self-scaled bound (SSB)}\label{sec:numerics-ssb} 

As explained in Section~\ref{sec:ssb}, to the best of our knowledge, no previous implementation of this bound is available; hence, Table~\ref{tab:newssb-all} reports only the behaviour of Algorithm~\ref{algo:ssb}.
For each instance, we report the total runtime and the resulting SSB value. The threshold was set to $\delta = 0.01$, which was found empirically to provide the fastest performance on the tested matrices.

\begin{table*}[!htbp]
\centering
\setlength{\tabcolsep}{4pt}
\renewcommand{\arraystretch}{1.05}

\caption{Results for the SSB obtained using Algorithm~\ref{algo:ssb} with $\epsilon = 10^{-6}$ and $\delta = 10^{-2}$. Running times are reported in seconds.}
\label{tab:newssb-all}

\begin{subtable}[t]{0.48\textwidth}
\centering
\footnotesize
\caption{Slack matrices of regular $n$-gons}
\label{tab:ssb-ngon}
\begingroup
\fontsize{8}{9}\selectfont
\begin{tabular}{lcc}
\toprule
Matrix & Time  & SSB \\
\midrule
$S_{5}$  & 6.48  & 5.000 \\
$S_{6}$  & 18.87  & 5.000 \\
$S_{7}$  & 59.76   & 5.542 \\
$S_{8}$  & 133   & 5.842 \\
$S_{9}$  & 1526    & 6.138 \\
$S_{10}$ & 928    & 6.410 \\
$S_{11}$ & 1869   & 6.669 \\
$S_{12}$ & 18902  & 6.918 \\
$S_{13}$ &  64108  & 7.133 \\
$S_{14}$ &  332085  & 7.314 \\
\bottomrule
\end{tabular}
\endgroup
\end{subtable}
\hfill
\begin{subtable}[t]{0.48\textwidth}
\centering
\footnotesize
\caption{Linear Euclidean distance matrices}
\label{tab:ssb-ledm}
\begingroup
\fontsize{8}{9}\selectfont
\begin{tabular}{lcc}
\toprule
Matrix & Time & SSB \\
\midrule
$M_{5}$ & 3.12   & 4.186 \\
$M_{6}$ & 8.02   & 4.593 \\
$M_{7}$ & 21.1   & 4.908 \\
$M_{8}$ & 230    & 5.251 \\
$M_{9}$ & 1007   & 5.506 \\
$M_{10}$ & 4609   & 5.763 \\
$M_{11}$ & 24219   & 6.004 \\
$M_{12}$ & 412601   & 6.219 \\
\bottomrule
\end{tabular}
\endgroup
\end{subtable}

\vspace{0.5em}

\begin{subtable}[t]{0.48\textwidth}
\centering
\footnotesize
\caption{Polyhedral slack matrices 
}
\label{tab:ssb-poly}
\begingroup
\fontsize{8}{9}\selectfont
\begin{tabular}{lcc}
\toprule
Matrix & Time & SSB \\
\midrule

Dodecahedron      & 5908 & 8.400  \\
24-cell           & 97200 & 10.667 \\
Icosidodecahedron & $\approx 1.9 \times 10^6$  & $10.243$ \\
Cuboctahedron     & 1491    & 7.800 \\
\bottomrule
\end{tabular}
\endgroup
\end{subtable}
\hfill
\begin{subtable}[t]{0.48\textwidth}
\centering
\footnotesize
\caption{Slack submatrices of the correlation polytope}
\label{tab:ssb-corr}
\begingroup
\fontsize{8}{9}\selectfont
\begin{tabular}{lcc}
\toprule
Matrix & Time & SSB \\
\midrule
$C_2$ & 0.777   & 3.667 \\
$C_3$ & 3.31    & 7.000 \\
$C_4$ & 469     & 12.791 \\
\bottomrule
\end{tabular}
\endgroup
\end{subtable}

\end{table*}

\paragraph{Discussion}

The results confirm that the SSB is more computationally demanding than the other bounds considered in this work.
Even on moderate-sized instances, the number of generated rank-one matrices grows quickly, and the total runtime becomes large. This is consistent with the nature of the separation problem, which is continuous, non-convex, and directly constrained by the entrywise values of $X$.
At the same time, the computed values show that SSB is often the strongest bound among the four approaches. On regular $n$-gon slack matrices and on LEDMs, the values obtained are significantly larger than those of the HSB and often match the best known lower bounds for the nonnegative rank. 
This confirms that exploiting the entrywise magnitude of the matrix, rather than only its support, can lead to stronger lower bounds. 

For the dodecahedron, $24$-cell and cuboctahedron, Algorithm~\ref{algo:ssb} yields numerical SSB estimates of $8.400$, $10.667$ and $7.800$ after $5908$, $97\,200$ and $1491$ seconds, respectively. For the icosidodecahedron, the numerical SSB estimate obtained after approximately 22 days is $10.243$. The computed lower bound $B>10$ yields $\rank_+\geq11$ (Table~\ref{tab:polytope-bounds}).

\subsection{New exact values and improved lower bounds for the benchmark matrices}\label{sec:numerics-new-results}

Beyond the computational comparisons reported in the previous sections, the different methods can be used as tools to discover several new exact values and improved lower bounds on the benchmark families considered in this paper. 
The complete and detailed tables are gathered in Appendix~\ref{app:tables}; we summarize here the main new results.

\begin{itemize}
    \item \textbf{LEDMs.}
    Our extended computations provide, to the best of our knowledge, the first numerical evaluations of the SSB for $M_{11}$ and $M_{12}$:
    
    \[
    \mathrm{SSB}(M_{11})\approx 6.004,
    \qquad
    \mathrm{SSB}(M_{12})\approx 6.219.
    \]
    
Both corrected lower bounds exceed $6$, recovering $\rank_+(M_{11})=\rank_+(M_{12})=7$.
    These ranks were already known from the geometric bound for $M_{10}$~\cite[Table~1]{Gillis_2012}, monotonicity under taking submatrices, and Hrube\v{s}' upper bounds~\cite{HRUBES2012457}.
    For both matrices, FSB, RCB, and HSB equal $3$, $6$, and $2$, respectively: SSB alone recovers the rank after taking the ceiling.

    \item \textbf{Regular $n$-gons.}
    For the rectangle covering bound, our non-convex formulation allows us to compute
    \[
    \mathrm{RCB}(S_{23})=8.
    \]
    To the best of our knowledge, this is a new exact value. As a consequence, the first value of $n$ for which $\mathrm{RCB}(S_n)$ is not known exactly is now $n=35$.

    \item \textbf{Slack matrices of selected polytopes.}
For the four polyhedral slack matrices considered here, we provide explicit exact nonnegative factorizations in the appendix, since such exact factorizations do not seem to be available in the literature. 
Combined with the known RCB values for the dodecahedron and the cuboctahedron, these factorizations close the gap and give the following exact extension complexity values
\[
\xc(\text{dodecahedron})=9
\qquad\text{and}\qquad
\xc(\text{cuboctahedron})=8.
\]
For the $24$-cell, we additionally computed $\mathrm{RCB}=12$ using the code of~\cite{fiorini2019}, which appears to be the most effective approach in our experiments (see Table~\ref{tab:rcb-extra-slack}). This matches the size of the exact factorization provided in the appendix, and therefore implies that the extension complexity of the $24$-cell is equal to $12$. To the best of our knowledge, this exact value was not previously available in the literature.

The numerical SSB estimates for the dodecahedron, $24$-cell and cuboctahedron are $8.400$, $10.667$ and $7.800$, respectively. For the icosidodecahedron, our SSB lower bound improves the previous nonnegative-rank lower bound $10$~\cite{dewez2021geometric} to $11$, leaving $11\leq\rank_+\leq12$.

    \item \textbf{UDISJ matrices.}
Using the code of~\cite{fiorini2019}, we were able to compute the rectangle covering bound for the matrices $U_n$ up to $n=6$:
\[
\mathrm{RCB}(U_1)=2, 
\mathrm{RCB}(U_2)=3, 
\mathrm{RCB}(U_3)=6, 
\mathrm{RCB}(U_4)=9, 
\mathrm{RCB}(U_5)=18, 
\mathrm{RCB}(U_6)=27.
\]
To the best of our knowledge, these are the largest exact values currently available for this family. Moreover, it means that for $n \leq 6$, the resulting sequence matches the upper bound on the RCB provided in~\cite{weltge2015sizes}, that is,
\[
\mathrm{RCB}(U_n)\le
\begin{cases}
(\sqrt{3})^n, & \text{if } n \text{ is even},\\[0.2em]
\dfrac{2}{\sqrt{3}}(\sqrt{3})^n, & \text{if } n \text{ is odd}.
\end{cases}
\]
If this holds for all $n \geq 7$, it would imply a stronger growth than the currently known general lower bound $\mathrm{RCB}(U_n)\geq (3/2)^n$~\cite{kaibel2015short}.

\item \textbf{Submatrices of the slack matrix of the correlation polytope.}
To the best of our knowledge, very little seems to be available in the literature specifically regarding numerical values of lower bounds for the matrices $C_n$.
The matrices $C_n$ are fully specified nonnegative completions of the UDISJ patterns $U_n$.
Hence, lower bounds on the FSB or RCB of the partial patterns also apply to $C_n$, although the exact values of these bounds may differ between the two families. Moreover, the authors of~\cite{Vandaele_2015} formulated the conjecture that $C_n$ has full nonnegative rank, that is,
\[
\rank_+(C_n)=2^n.
\]
In this context, the computations reported in Tables~\ref{tab:fsb-corr}, \ref{tab:rcb-corr}, \ref{tab:hsb-corr}, \ref{tab:ssb-corr}, and summarized in Table~\ref{tab:corrpolytope-bounds}, appear to provide the first published values for the four lower bounds considered in this work on this family. They imply in particular
\[
\rank_+(C_2)=4,\qquad \rank_+(C_3)\ge 7,
\]
which already follow from the ordinary rank, and
\[
\rank_+(C_4)\ge 13,\qquad
\rank_+(C_5)\ge 25,\qquad
\rank_+(C_6)\ge 34.
\]
For $n=4,5,6$, these improve substantially over the previously available generic lower bounds. Indeed, combining the ordinary-rank bound $1+\frac{n(n+1)}{2}$ with the bound $\lceil(3/2)^n\rceil$ inherited from the UDISJ pattern only gives 
\[
\rank_+(C_4)\geq 11,\qquad
\rank_+(C_5)\geq 16,\qquad
\rank_+(C_6)\geq 22.
\]

\end{itemize}

\section{Conclusion}\label{sec:conclusion}

In this paper, we revisited four classical lower bounds for the nonnegative rank: the fooling set bound (FSB), the rectangle covering bound (RCB), the hyperplane separation bound (HSB), and the self-scaled bound (SSB), and proposed exact non-convex optimization models and algorithms.
Our goal was to show that these bounds admit non-convex formulations that can be handled effectively by modern global solvers, such as Gurobi. 

For the HSB, we used an existing algorithmic framework from~\cite{fiorini2019} but reformulated the rank-one separation step as a non-convex optimization problem, thereby avoiding explicit MILP linearizations. 
Using the same principle, we introduced the first  practical algorithm for computing the SSB. 

The numerical experiments show that our non-convex approaches are not uniformly superior to the state of the art but provide a meaningful and often competitive alternative.
For the FSB and RCB, the proposed non-convex models are especially effective on dense matrices such as LEDMs and regular $n$-gon slack matrices. 
For the HSB, results must be mitigated: the MILP-based approaches remain faster on many instances. 
For the SSB, the proposed method is computationally  demanding, but it often provides the strongest lower bounds among the four approaches, and in some cases allowed us to close the gap with the best known upper bound. 
In view of the quality of the bounds obtained using SSB, our algorithm opens new perspectives for proving better lower bounds on the nonnegative rank of other matrices.

Because the practical performance of our non-convex formulations is closely tied to progress in global non-convex optimization, these approaches may scale significantly further in the future and extend the range of instances that can be solved to optimality.

It would be natural to investigate whether the same modeling philosophy can be extended to other lower bounds for the nonnegative rank, notably those based on nested-polytopes formulations~\cite{Gillis_2012,dewezphd} and refined rectangle-covering bounds such as the RRCB~\cite{oelze2014}.
Overall, the present work shows that non-convex global optimization provides not only effective computational tools but also a viable modeling framework for lower bounds on the nonnegative rank.

\revise{\paragraph*{Acknowledgments}

The authors thank Tilman Engel (University of Rostock) for identifying an error in the reported numerical results of the RCB of the Icosidodecahedron and for pointing out the omission of the geometric lower bounds of~\cite{Gillis_2012} from Table~\ref{tab:ngon-finite} for $n\geq24$ in an earlier version of this manuscript.}

{
\small
\microtypesetup{protrusion=false,expansion=false}
\titlespacing*{\section}{0pt}{3.5ex plus 1ex minus .2ex}{2.3ex plus .2ex}
\clubpenalty=150
\widowpenalty=150
 \bibliographystyle{abbrv}
 \bibliography{cas-refs}
}

\clearpage
\appendix

\setlist[itemize]{topsep=2pt,itemsep=1.5pt,parsep=0pt,partopsep=0pt,leftmargin=*}
\setlength{\abovedisplayskip}{5pt plus 2pt minus 2pt}
\setlength{\belowdisplayskip}{5pt plus 2pt minus 2pt}

\section{Up-to-date bounds for the nonnegative rank} \label{app:tables}

The tables summarize FSB, RCB, HSB, SSB, and the best known nonnegative-rank bounds for the benchmark families.
Numerical SSB values follow Section~\ref{sec:numerics-ssb}: completed estimates are validated to three decimal places by the computed enclosures. Integer rank bounds use the ceiling of the unrounded computed lower bound $B$, rather than the ceiling of the displayed estimate.

\subsection{Linear Euclidean Distance Matrices (LEDM)}\label{app:ledm}

Table~\ref{tab:ledm-finite} reports the results for the LEDMs.  
\begin{table}[H]
\centering
\footnotesize
\setlength{\tabcolsep}{6pt}
\renewcommand{\arraystretch}{1.08}
\begin{tabular}{c c c c c cc}
\toprule
& FSB
& RCB
& HSB
& SSB
& \multicolumn{2}{c}{$\rank_+(M_n)$} \\
\cmidrule(l){2-2} \cmidrule(l){3-3} \cmidrule(l){4-4} \cmidrule(l){5-5} \cmidrule(lr){6-7}
$n$ & & &  & & LB & UB \\
\midrule
$5$  & $3^{\mathrm{F13}}$ & $4^{\mathrm{D81}}$ & $2^{*}$ & $4.186^{*}$ & $5^{\mathrm{G12}}$ & $5^{\mathrm{triv}}$ \\
$6$  & $3^{\mathrm{F13}}$ & $4^{\mathrm{D81}}$ & $2^{*}$ & $4.593^{*}$ & $5^{\mathrm{G12}}$ & $5^{\mathrm{H12}}$ \\
$7$  & $3^{\mathrm{F13}}$ & $5^{\mathrm{D81}}$ & $2^{*}$ & $4.908^{*}$ & $6^{\mathrm{G12}}$ & $6^{\mathrm{H12}}$ \\
$8$  & $3^{\mathrm{F13}}$ & $5^{\mathrm{D81}}$ & $2^{*}$ & $5.251^{*}$ & $6^{\mathrm{G12}}$ & $6^{\mathrm{H12}}$ \\
$9$  & $3^{\mathrm{F13}}$ & $5^{\mathrm{D81}}$ & $2^{*}$ & $5.506^{*}$ & $6^{\mathrm{G12}}$ & $7^{\mathrm{H12}}$ \\
$10$ & $3^{\mathrm{F13}}$ & $5^{\mathrm{D81}}$ & $2^{*}$ & $5.763^{*}$ & $7^{\mathrm{G12}}$ & $7^{\mathrm{H12}}$ \\
$11$ & $3^{\mathrm{F13}}$ & $6^{\mathrm{D81}}$ & $2^{*}$ & $6.004^{*}$ & $7^{\mathrm{G12},\lceil\mathrm{SSB}\rceil}$ & $7^{\mathrm{H12}}$ \\
$12$ & $3^{\mathrm{F13}}$ & $6^{\mathrm{D81}}$ & $2^{*}$ & $6.219^{*}$ & $7^{\mathrm{G12},\lceil\mathrm{SSB}\rceil}$ & $7^{\mathrm{H12}}$ \\
\bottomrule
\end{tabular}
\caption{Summary of the values known for the LEDMs.}
\label{tab:ledm-finite}
\end{table}

\begin{flushleft}
\footnotesize
\emph{Sources.}
$^{*}$ = this work;
$\mathrm{F13}$ = \cite{fiorini2013combinatorial};
$\mathrm{D81}$ = \cite{de1981boolean};
$\mathrm{H12}$ = \cite{HRUBES2012457};
$\mathrm{G12}$ = \cite{Gillis_2012}; triv = trivial upper
bound.
In the $\rank_+$ columns, integer lower bounds derived from the SSB use the ceiling of the unrounded lower bound returned by Algorithm~\ref{algo:ssb}. For $n=11,12$, the superscript $\mathrm{G12}$ credits the previously known geometric lower bound, inherited from $M_{10}$ by monotonicity under taking submatrices; $\lceil\mathrm{SSB}\rceil$ records its independent numerical recovery using our SSB implementation.
\end{flushleft}

\normalsize

\begin{itemize}
    \item For the FSB, Fiorini et al.~\cite{fiorini2013combinatorial} proved that any nonnegative matrix with at most $s$ zero entries per row satisfies $\mathrm{FSB}(X)\le 2s+1$. Since each row of $M_n$ contains exactly one zero entry, this gives $\mathrm{FSB}(M_n)\leq 3$. Conversely, the three entries $(1,2)$, $(2,3)$, and $(3,1)$ always form a fooling set (see Figure~\ref{fig:ncfsLEDM}), hence $\mathrm{FSB}(M_n)\geq 3$. Therefore, $\mathrm{FSB}(M_n)=3$ for all $n\geq 3$.
    \item For the RCB, the exact value for the LEDM family follows from a result of De Caen~\cite{de1981boolean}. The lower bound is obtained via Sperner's theorem, while the upper bound comes from an explicit constructive covering, so that $\mathrm{RCB}(M_n)=\min\{r:\binom{r}{\lfloor r/2\rfloor}\ge n\}$.
    
    \item The HSB equals $2$ for all $n\geq2$ by Proposition~\ref{prop:hsb-ledm}, proved in Section~\ref{sec:numerics-hsb}.
    \item For the upper bound on $\rank_+(M_n)$, Hrube\v{s}~\cite{HRUBES2012457} proved that $\rank_+(M_{2n})\le \rank_+(M_n)+2$. Iterating this inequality $k$ times gives $\rank_+(M_n)\le \rank_+(M_{\lceil n/2^k\rceil})+2k$, and using the trivial bound $\rank_+(X)\le \min(m,n)$ for any nonnegative matrix $X\in\mathbb{R}_+^{m\times n}$, it gives $\rank_+(M_n)\leq \lceil n/2^k\rceil+2k$.
    Since this holds for every $k\in\mathbb{N}$, one obtains $\rank_+(M_n)\le \min_{k\in\mathbb{N}}(\lceil n/2^k\rceil+2k)$.    
    \item The SSB values for $M_{11}$ and $M_{12}$ are, to the best of our knowledge, reported here for the first time. After accounting for the separation optimality gap, they yield the sharp integer lower bound $7$ for both matrices. Their nonnegative ranks were already determined by combining the geometric lower bound of~\cite{Gillis_2012}, applied to the submatrix $M_{10}$, with the upper bounds of~\cite{HRUBES2012457}. Our computations show that the SSB also recovers these ranks, unlike the FSB, RCB, and HSB.
\end{itemize}

\subsection{Regular \texorpdfstring{$n$}{n}-gons}

Table~\ref{tab:ngon-finite} reports the results for the slack matrices of the regular $n$-gons.  

\begin{table}[ht!]
\centering
\scriptsize
\setlength{\tabcolsep}{12pt}
\renewcommand{\arraystretch}{1.08}
\begin{tabular}{c c cc c c cc}
\toprule
& FSB
& \multicolumn{2}{c}{RCB}
& HSB
& SSB
& \multicolumn{2}{c}{$\rank_+(S_n)$} \\
\cmidrule(l){2-2} \cmidrule(lr){3-4} \cmidrule(l){5-5} \cmidrule(l){6-6} \cmidrule(lr){7-8}
$n$ & & LB & UB &  & & LB & UB \\
\midrule
$5$        & $5$                  & $5^{\mathrm{FSB}}$   & $5^{\mathrm{triv}}$  & $3.090^{*}$ & $5^{*}$     & $5^{\mathrm{FSB}}$   & $5^{\mathrm{triv}}$ \\
$6$        & $4^{\mathrm{G20}}$   & $5^{\mathrm{V17}}$   & $5^{\mathrm{B86}}$   & $3.000^{*}$       & $5^{*}$     & $5^{\mathrm{RCB}}$   & $5^{\mathrm{V17}}$ \\
$7$        & $4^{\mathrm{G20}}$   & $6^{\mathrm{V17}}$   & $6^{\mathrm{B86}}$   & $3.115^{*}$ & $5.542^{*}$ & $6^{\mathrm{G12}}$   & $6^{\mathrm{F12}}$ \\
$8$        & $4^{\mathrm{G20}}$   & $6^{\mathrm{V17}}$   & $6^{\mathrm{B86}}$   & $3.314^{*}$ & $5.842^{*}$ & $6^{\mathrm{RCB}}$   & $6^{\mathrm{G00}}$ \\
$9$        & $4^{\mathrm{G20}}$   & $6^{\mathrm{V17}}$   & $6^{\mathrm{B86}}$   & $3.126^{*}$ & $6.138^{*}$ & $7^{\mathrm{O14}}$   & $7^{\mathrm{V17}}$ \\
$10$       & $4^{\mathrm{G20}}$   & $7^{\mathrm{V17}}$   & $7^{\mathrm{B86}}$   & $3.090^{*}$ & $6.410^{*}$ & $7^{\mathrm{G12}}$   & $7^{\mathrm{V17}}$ \\
$11$       & $4^{\mathrm{G20}}$   & $7^{\mathrm{V17}}$   & $7^{\mathrm{B86}}$   & $3.046^{*}$ & $6.669^{*}$ & $7^{\mathrm{G12}}$   & $7^{\mathrm{V17}}$ \\
$12$       & $4^{\mathrm{G20}}$   & $7^{\mathrm{V17}}$   & $7^{\mathrm{B86}}$   & $3.054^{*}$ & $6.918^{*}$ & $7^{\mathrm{G12}}$   & $7^{\mathrm{V17}}$ \\
$13$       & $4^{\mathrm{G20}}$   & $7^{\mathrm{V17}}$   & $7^{\mathrm{B86}}$   & $3.043^{*}$ & $7.133^{*}$ & $8^{\mathrm{O14}}$   & $8^{\mathrm{F12}}$ \\
$14$       & $4^{\mathrm{G20}}$   & $7^{\mathrm{V17}}$   & $7^{\mathrm{B86}}$   & $3.056^{*}$ & $7.314^{*}$ & $8^{\mathrm{O14}}$   & $8^{\mathrm{F12}}$ \\
$15$       & $4^{\mathrm{G20}}$   & $7^{\mathrm{V17}}$   & $7^{\mathrm{B86}}$   & $3.000^{*}$       & $-$         & $8^{\mathrm{G12}}$   & $8^{\mathrm{F12}}$ \\
$16$       & $4^{\mathrm{G20}}$   & $7^{\mathrm{V17}}$   & $7^{\mathrm{B86}}$   & $3.052^{*}$ & $-$         & $8^{\mathrm{G12}}$   & $8^{\mathrm{G00}}$ \\
$17$       & $4^{\mathrm{G20}}$   & $7^{\mathrm{V17}}$   & $7^{\mathrm{B86}}$   & $3.018^{*}$ & $-$         & $8^{\mathrm{G12}}$   & $9^{\mathrm{V17}}$ \\
$18$       & $4^{\mathrm{G20}}$   & $7^{\mathrm{V17}}$   & $7^{\mathrm{B86}}$   & $3.023^{*}$ & $-$         & $8^{\mathrm{G12}}$   & $9^{\mathrm{V17}}$ \\
$19$       & $4^{\mathrm{G20}}$   & $7^{\mathrm{V17}}$   & $7^{\mathrm{B86}}$   & $3.025^{*}$ & $-$         & $8^{\mathrm{G12}}$   & $9^{\mathrm{V17}}$ \\
$20$       & $4^{\mathrm{G20}}$   & $7^{\mathrm{V17}}$   & $7^{\mathrm{B86}}$   & $3.020^{*}$ & $-$         & $8^{\mathrm{G12}}$   & $9^{\mathrm{V17}}$ \\
$21$       & $4^{\mathrm{G20}}$   & $7^{\mathrm{V17}}$   & $7^{\mathrm{B86}}$   & $3.000^{*}$       & $-$         & $9^{\mathrm{G12}}$   & $9^{\mathrm{V17}}$ \\
$22$       & $4^{\mathrm{G20}}$   & $8^{\mathrm{D22}}$   & $8^{\mathrm{B86}}$   & $3.025^{*}$ & $-$         & $9^{\mathrm{G12}}$   & $9^{\mathrm{V17}}$ \\
$23$       & $4^{\mathrm{G20}}$   & $8^{*}$              & $8^{\mathrm{B86}}$   & $3.008^{*}$ & $-$         & $9^{\mathrm{G12}}$   & $9^{\mathrm{V17}}$ \\
$24$       & $4^{\mathrm{G20}}$   & $8^{\mathrm{V17}}$   & $8^{\mathrm{B86}}$   & $3.010^{*}$ & $-$         & $9^{\mathrm{G12}}$   & $9^{\mathrm{V17}}$ \\
$25$       & $4^{\mathrm{G20}}$   & $8^{\mathrm{V17}}$   & $8^{\mathrm{B86}}$   & $-$           & $-$         & $9^{\mathrm{G12}}$   & $10^{\mathrm{F12}}$ \\
$26$       & $4^{\mathrm{G20}}$   & $8^{\mathrm{V17}}$   & $8^{\mathrm{B86}}$   & $-$           & $-$         & $9^{\mathrm{G12}}$   & $10^{\mathrm{F12}}$ \\
$27$       & $4^{\mathrm{G20}}$   & $8^{\mathrm{V17}}$   & $8^{\mathrm{B86}}$   & $-$           & $-$         & $9^{\mathrm{G12}}$   & $10^{\mathrm{F12}}$ \\
$28$       & $4^{\mathrm{G20}}$   & $8^{\mathrm{V17}}$   & $8^{\mathrm{B86}}$   & $-$           & $-$         & $9^{\mathrm{G12}}$   & $10^{\mathrm{F12}}$ \\
$29$       & $4^{\mathrm{G20}}$   & $8^{\mathrm{V17}}$   & $8^{\mathrm{B86}}$   & $-$           & $-$         & $9^{\mathrm{G12}}$   & $10^{\mathrm{F12}}$ \\
$30$       & $4^{\mathrm{G20}}$   & $8^{\mathrm{V17}}$   & $8^{\mathrm{B86}}$   & $-$           & $-$         & $9^{\mathrm{G12}}$   & $10^{\mathrm{F12}}$ \\
$31$       & $4^{\mathrm{G20}}$   & $8^{\mathrm{V17}}$   & $8^{\mathrm{B86}}$   & $-$           & $-$         & $10^{\mathrm{G12}}$   & $10^{\mathrm{F12}}$ \\
$32$       & $4^{\mathrm{G20}}$   & $8^{\mathrm{V17}}$   & $8^{\mathrm{B86}}$   & $-$           & $-$         & $10^{\mathrm{G12}}$   & $10^{\mathrm{G00}}$ \\
$33$       & $4^{\mathrm{G20}}$   & $8^{\mathrm{V17}}$   & $8^{\mathrm{B86}}$   & $-$           & $-$         & $10^{\mathrm{G12}}$   & $11^{\mathrm{V17}}$ \\
$34$       & $4^{\mathrm{G20}}$   & $8^{\mathrm{V17}}$   & $8^{\mathrm{G16}}$   & $-$           & $-$         & $10^{\mathrm{G12}}$   & $11^{\mathrm{V17}}$ \\
$35$--$40$ & $4^{\mathrm{G20}}$   & $8^{\mathrm{V17}}$   & $9^{\mathrm{G16}}$   & $-$           & $-$         & $10^{\mathrm{G12}}$   & $11^{\mathrm{V17}}$ \\
$41$--$43$ & $4^{\mathrm{G20}}$   & $9^{\mathrm{V17}}$   & $9^{\mathrm{G16}}$   & $-$           & $-$         & $10^{\mathrm{G12}}$   & $11^{\mathrm{V17}}$ \\
$44$--$48$ & $4^{\mathrm{G20}}$   & $9^{\mathrm{V17}}$   & $9^{\mathrm{G16}}$   & $-$           & $-$         & $10^{\mathrm{G12}}$   & $11^{\mathrm{V17}}$ \\
$49$--$50$ & $4^{\mathrm{G20}}$   & $9^{\mathrm{V17}}$   & $9^{\mathrm{G16}}$   & $-$           & $-$         & $10^{\mathrm{G12}}$   & $12^{\mathrm{V17}}$ \\
$51$--$55$ & $4^{\mathrm{G20}}$   & $9^{\mathrm{V17}}$   & $9^{\mathrm{G16}}$   & $-$           & $-$         & $11^{\mathrm{G12}}$   & $12^{\mathrm{V17}}$ \\
$56$--$63$ & $4^{\mathrm{G20}}$   & $9^{\mathrm{V17}}$   & $10^{\mathrm{G16}}$  & $-$           & $-$         & $11^{\mathrm{G12}}$   & $12^{\mathrm{F12}}$ \\
$64$       & $4^{\mathrm{G20}}$   & $9^{\mathrm{V17}}$   & $10^{\mathrm{G16}}$  & $-$           & $-$         & $11^{\mathrm{G12}}$   & $12^{\mathrm{G00}}$ \\
$65$--$77$ & $4^{\mathrm{G20}}$   & $9^{\mathrm{V17}}$   & $10^{\mathrm{G16}}$  & $-$           & $-$         & $11^{\mathrm{G12}}$   & $13^{\mathrm{V17}}$ \\
$78$       & $4^{\mathrm{G20}}$   & $9^{\mathrm{V17}}$   & $10^{\mathrm{G16}}$  & $-$           & $-$         & $12^{\mathrm{G12}}$   & $13^{\mathrm{V17}}$ \\
$79$--$91$ & $4^{\mathrm{G20}}$   & $10^{\mathrm{V17}}$  & $10^{\mathrm{G16}}$  & $-$           & $-$         & $12^{\mathrm{G12}}$  & $13^{\mathrm{V17}}$ \\
\bottomrule
\end{tabular}
\caption{Summary of the values known for the slack matrices of the regular $n$-gons.\\ 
\emph{Sources:} 
$^{*}$ = this work;
$\mathrm{B86}$~=~\cite{barefoot1986biclique};
$\mathrm{G00}$~=~\cite{glineur2000computational};
$\mathrm{F12}$~=~\cite{Fiorini_2012};
$\mathrm{G12}$~=~\cite{Gillis_2012};
$\mathrm{O14}$~=~\cite{oelze2014};
$\mathrm{D22}$~=~\cite{dewezphd};
$\mathrm{V17}$~=~\cite{vandaele_2017};
$\mathrm{G16}$~=~\cite{GouchaGouveiaSilva2016};
$\mathrm{G20}$~=~\cite{gillis2020nonnegative};
$\mathrm{RCB}$ means that the lower bound on $\rank_+(S_n)$ is obtained from the RCB lower bound;
}
\label{tab:ngon-finite}
\end{table}

\begin{itemize}
    \item For the FSB, the general bound of Fiorini et al.~\cite{fiorini2013combinatorial} gives $\mathrm{FSB}(S_n)\leq 5$, since each row of $S_n$ contains exactly two zero entries. Goucha et al.~\cite{GouchaGouveiaSilva2016} state that the FSB equals $4$ for polygons, but they do not provide a proof.
    Such a proof is described in~\cite[Section~3.6.3.1]{gillis2020nonnegative}: the circulant zero pattern of $S_n$ excludes any $5\times 5$ submatrix with two zeros per row and column, and therefore implies $\mathrm{FSB}(S_n)=4$ for all $n\geq 6$.
    \item For the RCB, the pentagon is immediate since $\mathrm{RCB}(S_5)\ge \mathrm{FSB}(S_5)=5$. Beyond this case, and although small values were certainly computed informally within the community earlier, to the best of our knowledge the first published exact values of $\mathrm{RCB}(S_n)$ are those reported in~\cite{vandaele_2017} up to $n=13$. For larger $n$, the situation is different: the same paper provides a general lower bound on $\mathrm{RCB}(S_n)$, obtained by a refinement of Sperner's theorem, rather than exact values. On the upper-bound side, Barefoot et al.~\cite{barefoot1986biclique} reported Boolean-rank coverings numerically up to $n=33$, although, as emphasized in~\cite{GouchaGouveiaSilva2016}, no details were given to certify optimality; later, Goucha et al.~\cite{GouchaGouveiaSilva2016} extended this approach up to $n=91$ under the additional restriction that the rectangles have a prescribed homogeneous structure. For $n=22$, Dewez reports in his thesis~\cite{dewezphd} the exact value $\mathrm{RCB}(S_{22})=8$, obtained by solving the MILP formulation~\eqref{eq:rcb-dewez} in about $123000$ seconds, although this result does not seem to have appeared in a published paper. Finally, $\mathrm{RCB}(S_{23})=8$ is, to the best of our knowledge, a new exact value obtained here with our non-convex formulation; consequently, the first value of $n$ for which $\mathrm{RCB}(S_n)$ is not known exactly is now $n=35$.
    \item For the nonnegative rank, the cases $n=5$ and $n=6$ have long been known in the community, since the common lower bounds FSB and RCB already match the upper bounds there, although it is difficult to identify a first explicit reference.
    
    For the lower bounds, the first nontrivial case is $n=7$: using the geometric bound of Gillis and Glineur~\cite{Gillis_2012}, one obtains $\rank_+(S_7)\ge 6$, and the same argument also gives the tight lower bounds for $n=8,10,11,12,15,16$, $21\leq n\leq24$, and $n=31,32$, as well as the bounds $\rank_+(S_n)\geq 8$ for $n=17,18,19,20$; these values appear explicitly in~\cite{vandaele_2017}.
    The cases $n=9$ and $n=13$ were first settled at the values shown in the table using the refined rectangle covering bound of Oelze et al.~\cite{oelze2014}, as reported in~\cite{vandaele_2017}, while the case $n=14$ was later obtained by Dewez in his thesis~\cite{dewezphd} still using the refined rectangle covering bound.
    To obtain lower bounds for larger values of $n$ covered by Table~\ref{tab:ngon-finite}, we ran our own implementation of the geometric bound of Gillis and Glineur~\cite{Gillis_2012}. These computations yield $\rank_+(S_n)\geq 9$ for $24\leq n\leq30$, $\rank_+(S_n)\geq10$ for $31\leq n\leq50$, $\rank_+(S_n)\geq11$ for $51\leq n\leq77$, and $\rank_+(S_n)\geq12$ for $78\leq n\leq91$. Throughout this range, these bounds are stronger than the general RCB lower bound of~\cite{vandaele_2017}.

    For the upper bounds, Ben-Tal and Nemirovski~\cite{ben2001polyhedral} proved the bound $2\log_2(n)+4$ for powers of two, which was improved to $2\log_2(n)$ by Glineur~\cite{glineur2000computational}. For arbitrary $n$, Kaibel and Pashkovich~\cite{kaibel2011constructing,kaibel2013constructing} obtained the bound $2\lceil \log_2(n)\rceil+2$, later improved by Fiorini, Rothvo{\ss} and Tiwary~\cite{Fiorini_2012} to $2\lceil \log_2(n)\rceil$. Finally, Vandaele et al.~\cite{vandaele_2017} proved the sharper upper bound
$\rank_+(S_n)\le 2k-1$ for $2^{k-1}<n\le 2^{k-1}+2^{k-2}$ and $\rank_+(S_n)\le 2k$ for $2^{k-1}+2^{k-2}<n\le 2^k$, which has not been improved since and was conjectured to be exact.
\end{itemize}

\subsection{Slack matrices of selected polytopes}

Table~\ref{tab:polytope-bounds} reports the results for the slack matrices of four classical highly symmetric polytopes: the dodecahedron, the $24$-cell, the icosidodecahedron, and the cuboctahedron.

\begin{table}[H]
\centering
\scriptsize
\setlength{\tabcolsep}{8pt}
\renewcommand{\arraystretch}{1.08}
\begin{tabular}{l c cc cc cc cc}
\toprule
& FSB
& \multicolumn{2}{c}{RCB}
& \multicolumn{2}{c}{HSB}
& \multicolumn{2}{c}{SSB}
& \multicolumn{2}{c}{$\rank_+$} \\
\cmidrule(r){2-2} \cmidrule(lr){3-4} \cmidrule(lr){5-6} \cmidrule(lr){7-8} \cmidrule(lr){9-10}
Polytope &  & LB & UB & LB & UB & LB & UB & LB & UB \\
\midrule
Dodecahedron
& $6^{\mathrm{F19}}$
& $9^{\mathrm{F19}}$ & $9^{\mathrm{F19}}$
& $3.82^{\mathrm{F19}}$ & $3.82^{\mathrm{F19}}$
& $8.400^{*}$ & $8.400^{*}$
& $9^{\mathrm{RCB}}$ & $9^{\mathrm{V15*}}$ \\

24-cell
& $8^{\mathrm{F19}}$ 
& $12^{\mathrm{F19}}$ & $12^{\mathrm{F19}}$
& $5.6^{\mathrm{F19}}$ & $5.6^{\mathrm{F19}}$
& $10.667^{*}$ & $10.667^{*}$
& $12^{\mathrm{RCB}}$ & $12^{\mathrm{V15*}}$ \\

Icosidodecahedron
& $6^{\mathrm{F19}}$ 
& $8^{\mathrm{F19}}$ & $12^{\mathrm{F19}}$
& $3.708^{*}$ & $3.708^{*}$
& $10.243^{*}$ & $10.243^{*}$
& $11^{*}$ & $12^{\mathrm{V15*}}$ \\

Cuboctahedron
& $6^{\mathrm{F19}}$ 
& $8^{\mathrm{F19}}$ & $8^{\mathrm{F19}}$
& $4^{\mathrm{F19}}$ & $4^{\mathrm{F19}}$
& $7.800^{*}$ & $7.800^{*}$
& $8^{\mathrm{RCB}}$ & $8^{\mathrm{V15*}}$ \\
\bottomrule
\end{tabular}
\caption{Current values of FSB, RCB, HSB, SSB, and $\rank_+$ for selected polyhedral slack matrices. LB and UB denote lower and upper bounds, respectively. The SSB lower and upper bounds are rounded to the nearest $10^{-3}$ for display. Identical displayed bounds indicate agreement to three decimal places.\\
\emph{Sources:} 
$^{*}$ = this work;
$\mathrm{V15}$ = \cite{Vandaele_2015};
$\mathrm{V15}*$ = exact factorization presented in this work, obtained using the code of~\cite{Vandaele_2015};
$\mathrm{F19}$ = computed in this work using the C++ implementation of~\cite{fiorini2019};
$\mathrm{RCB}$ means that the lower bound on $\rank_+$ is obtained from the RCB lower bound.}
\label{tab:polytope-bounds}
\end{table}

\begin{itemize}
    \item For the slack matrices of the four polytopes considered here, some upper bounds on the extension complexity were known informally in the community, while others had only been obtained numerically, for instance with the heuristic code of \cite{Vandaele_2015}, but without published exact factorizations. Since the purpose of this appendix is to be as explicit as possible, we include below exact nonnegative factorizations, obtained by inspecting the numerical solutions returned by the heuristic of \cite{Vandaele_2015}, in order to provide analytic guarantees.
    \item For the dodecahedron, the value $\mathrm{RCB}=9$ has been known for some time in the community to match the extension complexity.
    The numerical SSB estimate is $8.400$. The computed lower bound $B>8$ independently yields the same integer lower bound $9$.
    An exact factorization of size $9$ is given by: 
{\scriptsize
\[
\setcounter{MaxMatrixCols}{24}
\renewcommand{\arraystretch}{0.9}
\setlength{\arraycolsep}{2.5pt}
W^\top = 
\begin{pmatrix}
1&0&0&1&0&0&1&c&0&0&b&1&0&0&b&c&0&0&0&0\\
0&1&0&0&b&0&0&0&1&b&0&0&c&0&0&0&c&1&1&0\\
1&0&c&0&0&0&1&0&1&0&0&0&0&c&0&b&b&0&1&0\\
0&1&b&0&c&b&1&0&0&0&c&1&0&0&0&0&0&0&1&0\\
0&1&0&1&0&c&0&b&0&0&0&1&b&0&0&0&0&1&0&c\\
1&0&0&1&0&0&0&0&1&c&0&0&0&b&c&0&0&1&0&b\\
0&0&0&0&1&0&0&0&0&1&1&0&0&0&1&0&0&0&0&0\\
0&0&1&0&0&1&0&0&0&0&0&0&0&1&0&0&0&0&0&1\\
0&0&0&0&0&0&0&1&0&0&0&0&1&0&0&1&1&0&0&0
\end{pmatrix}, 
\qquad
H=
\begin{pmatrix}
0&0&0&0&0&a&d&0&d&a&0&0\\
0&0&0&d&a&0&0&0&0&0&a&d\\
0&a&0&0&d&d&0&a&0&0&0&0\\
d&0&0&0&0&0&0&d&a&0&0&a\\
a&0&a&0&0&0&0&0&0&d&d&0\\
0&d&d&a&0&0&a&0&0&0&0&0\\
0&0&0&0&a&a&0&0&0&a&a&0\\
0&0&0&a&0&0&a&0&a&0&0&a\\
a&a&a&0&0&0&0&a&0&0&0&0
\end{pmatrix}. \] }
 where 
$\phi=\frac{1+\sqrt5}{2}$, 
$a=\frac{2}{\phi^2}$, 
$b=\frac{1}{\phi}$, 
$c=\phi$, and 
$d=2b$.

    \item For the $24$-cell, Algorithm~\ref{algo:ssb} yields a numerical SSB estimate of $10.667$ after $97\,200$ seconds. On the other hand, using the code of Fiorini et al.~\cite{fiorini2019}, we computed that $\mathrm{RCB}=12$, a value for which we did not find an explicit trace in the literature.
    Below, we provide an exact nonnegative factorization of size $12$, which certifies that the extension complexity of the $24$-cell is $12$:

{\scriptsize
\[
\setcounter{MaxMatrixCols}{24}
\renewcommand{\arraystretch}{0.9}
\setlength{\arraycolsep}{2.5pt}
W=
\begin{pmatrix}
0&0&0&1&0&0&0&1&1&1&0&0\\
1&0&0&1&0&0&0&0&1&0&1&0\\
0&0&1&1&0&0&0&0&0&1&1&0\\
0&0&0&1&0&1&0&0&0&1&1&0\\
0&0&0&1&0&0&1&0&1&0&1&0\\
1&0&0&0&0&0&0&1&1&0&0&1\\
0&0&1&0&0&0&0&1&0&1&0&1\\
0&0&0&0&0&1&0&1&0&1&0&1\\
0&0&0&0&0&0&1&1&1&0&0&1\\
0&1&0&1&0&0&0&0&1&1&0&0\\
1&0&1&0&0&0&0&0&0&0&1&1\\
1&0&0&0&0&1&0&0&0&0&1&1\\
0&0&1&0&0&0&1&0&0&0&1&1\\
0&0&0&0&0&1&1&0&0&0&1&1\\
0&0&0&0&1&0&0&1&1&1&0&0\\
1&1&0&0&0&0&0&0&1&0&0&1\\
0&1&1&0&0&0&0&0&0&1&0&1\\
0&1&0&0&0&1&0&0&0&1&0&1\\
0&1&0&0&0&0&1&0&1&0&0&1\\
1&0&0&0&1&0&0&0&1&0&1&0\\
0&0&1&0&1&0&0&0&0&1&1&0\\
0&0&0&0&1&1&0&0&0&1&1&0\\
0&0&0&0&1&0&1&0&1&0&1&0\\
0&1&0&0&1&0&0&0&1&1&0&0
\end{pmatrix},
\qquad
H=
\begin{pmatrix}
0&0&0&0&0&0&0&1&1&1&1&0&0&0&0&0&0&0&0&1&1&1&1&2\\
0&1&1&0&0&0&2&1&1&0&0&0&0&0&0&1&1&0&0&0&0&1&1&0\\
0&0&1&0&1&0&0&1&0&1&0&0&0&0&1&0&1&2&0&1&0&1&0&0\\
0&0&0&0&0&0&0&0&0&0&0&0&0&1&1&1&1&0&2&1&1&1&1&0\\
2&1&1&1&1&0&0&1&1&1&1&0&0&0&0&0&0&0&0&0&0&0&0&0\\
0&1&0&1&0&0&0&0&1&0&1&2&0&1&0&1&0&0&0&0&1&0&1&0\\
0&1&1&1&1&2&0&0&0&0&0&0&0&1&1&1&1&0&0&0&0&0&0&0\\
0&0&0&1&1&0&0&0&0&1&1&0&2&1&1&0&0&0&0&1&1&0&0&0\\
0&0&0&0&0&0&0&0&0&0&0&1&0&0&0&0&0&1&0&0&0&0&0&0\\
0&0&0&0&0&1&0&0&0&0&0&0&0&0&0&0&0&0&0&0&0&0&0&1\\
0&0&0&0&0&0&1&0&0&0&0&0&1&0&0&0&0&0&0&0&0&0&0&0\\
1&0&0&0&0&0&0&0&0&0&0&0&0&0&0&0&0&0&1&0&0&0&0&0
\end{pmatrix}.
\]}
    \item For the icosidodecahedron, the code of Fiorini et al.~\cite{fiorini2019} gives $\revise{8}\leq\mathrm{RCB}\leq\revise{12}$.The numerical SSB estimate is $10.243$. The computed lower bound $B>10$ yields $\rank_+\geq11$, improving the previous lower bound $10$~\cite{dewez2021geometric}.
    Together with the explicit exact nonnegative factorization of size $12$ below, this gives $11\leq\rank_+\leq12$:

{\scriptsize
\setcounter{MaxMatrixCols}{32}
\renewcommand{\arraystretch}{0.88}
\setlength{\arraycolsep}{2.2pt}

\[
W^\top=
\begin{pmatrix}
0&a&0&a&0&b&a&b&b&1&0&0&0&1&a&1&b&0&0&0&0&1&d&0&0&0&0&0&0&0\\
0&1&0&1&0&a&0&0&0&b&b&a&1&0&0&b&a&0&1&a&b&0&a&0&0&a&0&d&0&0\\
0&b&a&0&0&b&0&a&0&0&0&0&0&0&0&a&0&0&b&0&a&0&0&a&0&b&0&a&0&1\\
0&1&0&0&0&d&a&1&0&a&a&0&0&b&0&1&0&a&1&b&1&a&b&1&b&d&0&a&a&e\\
b&0&1&0&a&0&0&0&b&0&0&a&0&0&b&0&a&a&0&0&0&a&0&0&0&0&b&0&0&a\\
a&0&0&0&0&0&0&0&0&0&1&b&a&0&0&0&0&1&a&d&1&0&0&b&1&b&b&0&a&0\\
1&0&0&0&d&0&1&a&a&0&0&0&0&b&1&0&0&b&0&a&0&b&a&a&b&0&a&0&1&0\\
0&0&b&b&0&0&0&0&a&a&a&b&b&0&0&0&b&0&0&0&0&0&0&0&0&0&a&a&0&0\\
a&a&f&0&0&a&0&f&a&a&a&f&0&a&a&b&f&b&a&b&b&b&b&f&a&a&a&0&0&f\\
0&0&0&0&a&a&b&b&0&0&0&0&0&a&0&0&0&0&0&0&0&0&0&b&a&a&0&0&b&b\\
a&a&0&a&b&0&a&0&0&0&0&0&a&0&a&0&0&0&a&0&0&0&0&0&0&0&0&b&a&0\\
1&0&e&a&a&0&0&0&d&b&b&1&a&a&1&a&1&1&0&b&a&1&b&0&a&0&d&0&0&0
\end{pmatrix}
\]

\[
H=
\begin{pmatrix}
0&b&1&0&0&0&0&1&c&0&0&0&c&1&0&0&0&0&0&0&1&b&1&0&0&0&c&0&0&0&1&c\\
1&0&1&1&1&0&0&0&0&1&0&0&0&1&1&0&0&0&0&0&0&0&1&1&0&0&0&1&0&0&0&0\\
0&0&0&a&0&0&0&0&1&0&1&0&0&1&1&0&0&0&1&0&1&0&1&1&0&0&0&0&1&0&0&1\\
0&0&0&0&0&0&0&0&0&1&0&0&0&0&0&0&1&1&0&0&0&1&0&0&0&0&0&1&0&0&0&0\\
0&0&1&0&1&0&1&1&0&0&0&0&1&0&0&1&0&0&0&a&0&0&0&0&1&0&1&0&0&1&1&0\\
0&0&0&0&1&b&1&0&0&0&c&0&0&0&1&c&0&b&1&0&0&0&0&1&c&0&0&0&c&1&0&0\\
0&0&0&0&0&0&1&1&0&0&0&1&0&0&0&0&1&0&1&1&1&0&0&0&0&1&0&0&0&1&1&0\\
1&1&1&c&1&1&1&1&1&1&1&1&1&1&1&1&0&0&0&0&0&0&0&0&0&0&0&0&0&0&0&0\\
0&0&0&1&0&0&0&0&0&0&0&0&0&0&0&0&0&0&0&1&0&0&0&0&0&0&0&0&0&0&0&0\\
0&0&0&0&0&0&0&0&0&0&0&0&0&0&0&0&1&1&1&c&1&1&1&1&1&1&1&1&1&1&1&1\\
0&1&0&0&0&1&0&0&1&0&1&0&1&0&0&1&0&1&0&0&0&1&0&0&1&0&1&0&1&0&0&1\\
1&1&0&0&0&1&0&0&0&0&0&1&0&0&0&0&0&0&0&0&0&0&0&0&0&1&0&0&0&0&0&0
\end{pmatrix}.
\]
}
where 
$\phi=\frac{1+\sqrt5}{2}$, 
$a=\frac{1}{\phi^2}$, 
$b=\frac{1}{\phi}$, 
$c=\phi$, 
$d=2b$, 
$e=1+a$, and 
$f=b-a$.

    \item For the cuboctahedron, Dewez et al.~\cite{dewez2021geometric} report that $\mathrm{RCB}=8$.
    The numerical SSB estimate is $7.800$. The computed lower bound $B>7$ independently yields $\rank_+\geq8$, matching the size of the exact nonnegative factorization given below. Thus, the extension complexity of the cuboctahedron is $8$:
{\scriptsize
\[
\setcounter{MaxMatrixCols}{14}
\renewcommand{\arraystretch}{0.9}
\setlength{\arraycolsep}{2.5pt}
W^\top = 
\begin{pmatrix}
1&0&2&0&1&1&0&0&1&0&0&0\\
1&0&1&1&1&0&1&0&0&1&0&0\\
0&0&0&1&0&0&1&1&0&2&0&1\\
0&0&0&0&1&0&1&0&1&0&2&1\\
0&1&0&0&0&1&0&1&1&0&1&1\\
1&2&0&1&0&1&0&1&0&0&0&0\\
1&1&0&1&1&0&1&0&0&0&1&0\\
0&0&1&0&0&1&0&1&1&1&0&1
\end{pmatrix}, 
\qquad
H=
\begin{pmatrix}
1&0&0&0&0&1&0&0&0&0&1&0&0&1\\
0&0&1&0&0&0&1&0&0&0&0&0&1&0\\
0&1&0&1&0&0&0&1&1&0&0&0&0&0\\
1&0&0&1&0&0&0&0&0&0&0&1&1&0\\
0&0&0&0&1&0&0&1&0&0&0&0&0&1\\
0&1&0&0&0&1&1&0&0&1&0&0&0&0\\
0&0&1&0&0&0&0&0&1&0&1&0&0&0\\
0&0&0&0&1&0&0&0&0&1&0&1&0&0
\end{pmatrix}.
\]}

\end{itemize}

\subsection{Slack submatrices associated with the correlation polytope}

Table~\ref{tab:corrpolytope-bounds} reports the results for $C_n$. 

\begin{table}[ht!]
\centering
\footnotesize
\setlength{\tabcolsep}{5pt}
\renewcommand{\arraystretch}{1.08}
\begin{tabular}{c cc cc c c cc}
\toprule
& \multicolumn{2}{c}{FSB}
& \multicolumn{2}{c}{RCB}
& HSB
& SSB
& \multicolumn{2}{c}{$\rank_+(C_n)$} \\
\cmidrule(lr){2-3} \cmidrule(lr){4-5} \cmidrule(l){6-6} \cmidrule(l){7-7} \cmidrule(lr){8-9}
$n$ & LB & UB & LB & UB &  &  & LB & UB \\
\midrule
$2$ & $3^{\mathrm{F19}}$  & $3^{\mathrm{F19}}$  & $3^{\mathrm{F19}}$  & $3^{\mathrm{F19}}$  & $3.667^{\mathrm{F19}}$  & $3.667^{*}$    & $4^{\mathrm{rank}}$ & $4^{\mathrm{triv}}$ \\
$3$ & $7^{\mathrm{F19}}$  & $7^{\mathrm{F19}}$  & $7^{\mathrm{F19}}$  & $7^{\mathrm{F19}}$  & $2.500^{\mathrm{F19}}$     & $7.000^{*}$    & $7^{\mathrm{RCB}}$     & $8^{\mathrm{triv}}$ \\
$4$ & $10^{\mathrm{F19}}$ & $10^{\mathrm{F19}}$ & $13^{\mathrm{F19}}$ & $13^{\mathrm{F19}}$ & $2.032^{\mathrm{F19}}$  & $12.791^{*}$   & $13^{\mathrm{RCB}}$    & $16^{\mathrm{triv}}$ \\
$5$ & $16^{\mathrm{F19}}$ & $16^{\mathrm{F19}}$ & $25^{\mathrm{F19}}$ & $25^{\mathrm{F19}}$ & $2.100^{\mathrm{F19}}$  & $-$            & $25^{\mathrm{RCB}}$    & $32^{\mathrm{triv}}$ \\
$6$ & $19^{\mathrm{F19}}$ & $31^{\mathrm{F19}}$ & $34^{\mathrm{F19}}$ & $39^{\mathrm{F19}}$ & $-$                      & $-$            & $34^{\mathrm{RCB}}$    & $64^{\mathrm{triv}}$ \\
$n$ & $-$                 & $-$                 & $\left\lceil(3/2)^n\right\rceil^{\mathrm{K15}}$ & $2^n{}^{\mathrm{triv}}$ & $-$ & $-$ & $\max\!\left\{\rank(C_n), \left\lceil(3/2)^n\right\rceil \right\}^{\mathrm{rank},\mathrm{RCB}}$ & $2^n{}^{\mathrm{triv}}$ \\
\bottomrule
\end{tabular}
\caption{Current values of FSB, RCB, HSB, SSB, and $\rank_+$ for slack submatrices associated with the correlation polytope. For each quantity, LB and UB denote the best currently known lower and upper bounds, respectively.\\ 
\emph{Sources:} 
$^{*}$ = this work;
$\mathrm{K15}$ = \cite{kaibel2015short};
$\mathrm{F19}$ = computed in this work using the C++ implementation of~\cite{fiorini2019};
$\mathrm{RCB}$ means that the lower bound on $\rank_+(C_n)$ is obtained from the RCB lower bound; $\mathrm{rank}$ = ordinary-rank lower bound;
$\mathrm{triv}$ = trivial upper bound $2^n$. 
}
\label{tab:corrpolytope-bounds}
\end{table}

Since $C_n$ is a nonnegative completion of $U_n$, the bound of~\cite{kaibel2015short} gives $\mathrm{RCB}(C_n)\geq(3/2)^n$.
For $C_4$, the unrounded SSB lower bound $B>12$ also yields $\rank_+(C_4)\geq13$.

\end{document}